\documentclass[11pt]{article}

\usepackage{graphicx}
\usepackage{amsmath,amssymb,amsthm}
\usepackage{mathtools}
\usepackage[margin=3cm]{geometry}
\usepackage{booktabs}
\usepackage{float}
\usepackage{enumitem}
\floatstyle{ruled}
\newfloat{algorithm}{htbp}{loa}
\floatname{algorithm}{Algorithm}
\usepackage{tikz-cd}
\usepackage{mathrsfs}
\usepackage[expansion=false]{microtype}
\usepackage[hidelinks]{hyperref}

\newcommand{\R}{\mathbb{R}}

\newcommand{\diag}{\operatorname{diag}}
\newcommand{\tr}{\operatorname{tr}}
\DeclareMathOperator{\spann}{span}
\DeclareMathOperator{\ad}{ad}
\DeclareMathOperator{\Ad}{Ad}
\DeclareMathOperator{\cond}{cond}
\DeclareMathOperator{\inj}{inj}
\DeclareMathOperator{\supp}{supp}

\DeclareMathOperator{\dist}{dist}
\DeclareMathOperator{\grad}{grad}
\DeclareMathOperator{\Hess}{Hess}

\newcommand{\Sym}{\operatorname{Sym}}
\newcommand{\Spd}{\mathcal{P}}
\newcommand{\g}{\mathfrak{g}}
\newcommand{\dG}{d_G}
\newcommand{\dAI}{d_{\mathrm{AI}}}
\newcommand{\dL}{d_{L}}
\newcommand{\Lam}{\Lambda}
\newcommand{\Th}{\Theta}
\newcommand{\Nat}[1]{\mathcal{N}_{#1}}
\newcommand{\Hk}{H_K}

\newcommand{\ip}[2]{\left\langle #1, #2 \right\rangle}
\newcommand{\Ceq}{C_{\mathrm{eq}}}
\newcommand{\cqu}{c_{\mathrm{qu}}}
\newcommand{\Llog}{\mathscr{L}}
\newcommand{\Eexp}{\mathscr{E}}

\newtheorem{theorem}{Theorem}[section]
\newtheorem{corollary}[theorem]{Corollary}
\newtheorem{lemma}[theorem]{Lemma}
\newtheorem{proposition}[theorem]{Proposition}
\newtheorem{definition}[theorem]{Definition}
\newtheorem{assumption}[theorem]{Assumption}

\newtheorem{example}{Example}[section]
\newtheorem{remark}[theorem]{Remark}

\usepackage{xcolor}

\title{Multilevel kernel approximation of log-Euclidean manifold-valued functions}

\author{
  Nir Sharon\thanks{Email: \texttt{nsharon@tauex.tau.ac.il}} \\
  \small School of Mathematical Sciences \\
  \small Tel Aviv University \\
  \small Tel Aviv 6997801, Israel
  \and
  Holger Wendland\thanks{Email: \texttt{holger.wendland@uni-bayreuth.de}} \\
  \small Department of Mathematics \\
  \small University of Bayreuth \\
  \small 95440 Bayreuth, Germany
}

\date{}

\begin{document}
\maketitle

\begin{abstract}
We study kernel approximation and multilevel residual correction for functions with values in a manifold that admits a global logarithm onto a finite-dimensional Hilbert space, such as positive definite matrices, exponential Lie groups, and Hadamard manifolds. Approximating in logarithmic coordinates keeps the approximant on the manifold and transfers Hilbert-space error estimates verbatim. To express these errors in the natural geometry of the target, we prove metric comparisons. In particular, we obtain a sharp, scale-invariant bound between the log-Euclidean and affine-invariant distances on positive definite matrices, with a constant that depends only on condition numbers. For multilevel approximation, compactly supported kernels scaled to the fill distance give algebraic convergence with sparse, uniformly well-conditioned systems. These guarantees require appropriate sample locations and approximation scales. Compact groups are treated through a principal logarithm, where a bi-invariant metric gives error bounds with constant one. Finally, the log-Euclidean and Karcher means differ by $\mathcal O(r^3)$ near the identity and by $\mathcal O(r^2)$ around a general point, where $r$ is the logarithmic data radius and the second bound is sharp.
\end{abstract}

\section{Introduction} \label{sec:intro}

Approximating a function from scattered samples is a basic task of numerical analysis, and the classical theory deals with functions whose values lie in a vector space. In a growing number of applications, the values lie instead in a nonlinear set: rotations and rigid motions in robotics and computer vision, deformation and diffusion tensors in mechanics and medical imaging, and covariance matrices in statistics take values in matrix Lie groups or in the cone $\Spd_n$ of symmetric positive definite matrices. Treating such data componentwise in the ambient matrix space is not an option: a kernel interpolant of positive definite matrices need not be positive definite, and even where it is, one observes the well-documented swelling of determinants \cite{Arsigny2006, Arsigny2007}; a linear average of rotations is not a rotation. The methods proposed for manifold-valued data fall into three groups, see the survey \cite{Zimmermann2021}: intrinsic schemes built from Riemannian (Karcher) means, schemes that work in a tangent space or in a chart, and projection schemes that approximate in an ambient space and project back onto the manifold, as in \cite{FischerWendland2026}. This paper belongs to the second group.

We consider the same construction for a smooth manifold $G$ with a diffeomorphism $\log:G\to H$ onto a finite-dimensional Hilbert space. Such pullback geometries have recently been studied in \cite{BissonPennec2025}; examples include positive definite matrices, exponential Lie groups, and Hadamard manifolds. The identity
\[
 d_G(g,h):=\|\log g-\log h\|_H
\]
transfers approximation estimates from $H$ to $G$. Our main geometric result makes this useful for natural distances. For positive definite matrices, we prove a sharp comparison between the log-Euclidean and affine-invariant distances. This comparison converts approximation estimates in logarithmic coordinates into error bounds in the natural geometry. The distortion depends only on the matrices’ condition numbers and is unchanged by scalar rescaling. We also give comparisons for invariant metrics on Lie groups and for Riemannian metrics using Rauch's theorem.

The approximation method combines the kernel construction of \cite{WendlandSPD} with the multilevel residual correction of \cite{sharon2022multiscale}. We further analyze the bias-reducing effect of multilevel residual correction for scalar- and manifold-valued quasi-interpolation in~\cite{AbasSharon2025}. Here, we prove that the transported scheme is conjugate, level by level, to residual correction in $H$. We also formalize and verify selected algebraic results in Lean~4, as documented in the accompanying code. Since a fixed kernel with nested sites reduces to interpolation on the finest level, an elementary observation recorded in Section~\ref{subsec:telescope}, we use compactly supported kernels scaled to the fill distance, for which the scalar convergence theory \cite{GollwitzerKempfWendland2026, Wendland2010} extends componentwise to $H$-valued functions. We state this transfer with its contraction factor and computational cost. These convergence results are essentially inherited from scalar approximation. However, we focus mainly on geometric comparisons, their use for natural-metric errors, and comparisons with intrinsic averages. 

Compact groups require a local logarithm, but a bi-invariant metric still gives a global one-sided error bound with constant one. We further compare the log-Euclidean mean with the Karcher mean: their distance is $\mathcal O(r^3)$ for logarithms of size $r$ near zero, and $\mathcal O(r^2)$ for data clustered around a general point. An explicit positive definite example attains the theoretical second order. 

The paper is organized as follows. Sections~\ref{sec:targets}--\ref{sec:transference} develop the geometry, Sections~\ref{sec:singlescale}--\ref{sec:computation} give the approximation and cost estimates, and Sections~\ref{sec:local}--\ref{sec:numerics} treat compact groups, intrinsic means and numerical examples. Proofs of the metric comparisons and multilevel bound are in the appendices.

\section{Log-Hilbert targets} \label{sec:targets}

Throughout, $H$ is a real Hilbert space of finite dimension $m$ with inner product $\ip{\cdot}{\cdot}_H$ and norm $\|\cdot\|_H$.

\begin{definition}
\label{def:loghilbert}
A \emph{log-Hilbert target} with model space $H$ is a smooth manifold $G$ together with a diffeomorphism $\log:G\to H$. We write $\exp:=\log^{-1}$ and $e:=\exp(0)$. If $G$ is a Lie group with Lie algebra $\g$, we always take $H=\g$ with some inner product and $\log$ the inverse of the Lie-theoretic exponential map $\exp_G:\g\to G$; the definition then requires $\exp_G$ to be a diffeomorphism, and such groups are called \emph{exponential}.
\end{definition}
This is also the setting of~\cite{BissonPennec2025}, where $G$ is called a log-Euclidean Lie group with respect to the transported addition. The chart described in the above definition is part of the target structure. Its smoothness identifies intrinsic and extrinsic regularity on compact subsets, as we will see in Lemma~\ref{lem:moser}. Nevertheless, it need not respect an existing multiplication. Next, we provide three examples. For positive definite matrices of the first example, the transported addition is the log-Euclidean operation of \cite{Arsigny2007}. In particular, for a non-Abelian Lie group, it differs from group multiplication.

\begin{example}[Positive definite matrices]
\label{ex:spd}
The cone $\Spd_n$ of symmetric positive definite $n\times n$ matrices is a log-Hilbert target with $H=\Sym_n$, the symmetric matrices with the Frobenius inner product $\ip{X}{Y}=\tr(XY)$, and $\log$ the principal matrix logarithm; $\exp:\Sym_n\to\Spd_n$ is an analytic diffeomorphism \cite[Ch.~11]{Higham2008}. This is the log-Euclidean structure of \cite{Arsigny2007} and the setting of \cite{WendlandSPD}.
\end{example}

\begin{example}[Exponential Lie groups]
\label{ex:lie}
Every connected, simply connected nilpotent Lie group is exponential \cite[Thm.~1.2.1]{CorwinGreenleaf1990}, including the Heisenberg groups and unipotent upper triangular groups $\mathrm{UT}_n$. There $\log(I+N)=\sum_{k=1}^{n-1}(-1)^{k+1}N^k/k$. More generally, a connected, simply connected solvable Lie group is exponential exactly when no $\ad_X$ has a nonzero purely imaginary eigenvalue \cite{Dixmier1957,HilgertNeeb2012}. Upper triangular matrices with positive diagonal and their connected Lie subgroups provide further examples.
\end{example}

\begin{example}[Hadamard manifolds]
\label{ex:hadamard}
For a complete, simply connected Riemannian manifold of nonpositive sectional curvature, the Cartan--Hadamard theorem makes $\exp_p:T_pM\to M$ a diffeomorphism. Thus $H=T_pM$ and $\log=\exp_p^{-1}$ define a log-Hilbert target, depending on $p$. For $\Spd_n$ with its affine-invariant metric, $\exp_I$ is the matrix exponential \cite[Ch.~6]{Bhatia2007}. General CAT(0) spaces, such as metric trees, need not have a linear tangent space and require other constructions.
\end{example}

Finite products of log-Hilbert targets are log-Hilbert targets. Compact connected Lie groups of positive dimension are not, since $\exp$ is not injective; neither are $\mathrm{GL}^+(n)$ and $\mathrm{SL}(n)$ for $n\ge2$, where $\exp$ is not surjective, nor compact manifolds such as spheres and Grassmannians. Section~\ref{sec:local} treats compact groups by means of a principal logarithm.

We now transport the structure of $H$ to $G$. For $g,h\in G$ and $\lambda\in\R$ set
\begin{align*}
    g \oplus h &:= \exp\big(\log g + \log h\big), &
    \lambda \odot g &:= \exp\big(\lambda\,\log g\big), &
    g\ominus h &:=\exp(\log g-\log h),\\
    \ip{g}{h}_G &:= \ip{\log g}{\log h}_H, &
    \|g\|_G &:= \|\log g\|_H, &
    \dG(g,h) &:= \|\log g - \log h\|_H .
\end{align*}

\begin{proposition}
\label{prop:hilbert}
With these operations, $G$ is a real Hilbert space, $\exp:(H,+,\cdot)\to (G,\oplus,\odot)$ is a linear isometric isomorphism, and $(G,\dG)$ is a complete metric space.
\end{proposition}

\begin{proof}
For $\Spd_n$, this is \cite[Thm.~2.2]{WendlandSPD}; the argument, transport of structure through the bijection $\exp$, is the same in general. See also \cite{BissonPennec2025}.
\end{proof}

The Baker--Campbell--Hausdorff formula gives $\log(gh)=\log g+\log h+\tfrac12[\log g,\log h]+\cdots$ locally, so $\oplus$ generally differs from Lie-group multiplication. The relevant approximation error is often measured instead by an invariant Riemannian distance, the affine-invariant distance on $\Spd_n$, or an ambient matrix norm. We next compare these with $d_G$.

\section{The transference principle and the natural metrics} \label{sec:transference}

\subsection{Function spaces and the transference identity}

Let $D\subset\R^d$ be a bounded domain with Lipschitz boundary,
and fix an orthonormal basis $w_1,\dots,w_m$ of $H$.
We write $\mathcal N(D)$ for a scalar function space on $D$.
The spaces considered here are $C^k(\overline D)$, $L_q(D)$,
the Sobolev spaces $W^\tau_q(D)$, including
$H^\sigma(D)=W^\sigma_2(D)$, and reproducing kernel Hilbert spaces.

The corresponding space $\mathcal N(D;H)$ consists of
$H$-valued functions
\[
u=\sum_{\alpha=1}^m u_\alpha w_\alpha,
\qquad u_\alpha\in\mathcal N(D).
\]
When $\mathcal N(D)$ is a Hilbert space, we use the norm
\[
\|u\|_{\mathcal N(D;H)}^2
=\sum_{\alpha=1}^m\|u_\alpha\|_{\mathcal N(D)}^2.
\]
For the other spaces, we use the standard vector-valued norm, replacing the scalar absolute value with $\|\cdot\|_H$ in its definition. Since $H$ is finite-dimensional, these norms are equivalent to norms constructed from the scalar component norms. For $H=\Sym_n$ with the Frobenius inner product, this gives the componentwise function spaces used in \cite[Sec.~4]{WendlandSPD}.

The corresponding spaces of $G$-valued functions are defined by pullback,
\[
  \mathcal N(D;G) := \{\, F:D\to G \ :\ \log\circ F\in \mathcal N(D;H)\,\}, \qquad
  \|F\|_{\mathcal N(D;G)} := \|\log\circ F\|_{\mathcal N(D;H)},
\]
and the discrepancy between two $G$-valued functions is measured through $\ominus$ pointwise,
\[
  \|F\ominus \tilde F\|_{\mathcal N(D;G)}
  := \|\log\circ F - \log\circ \tilde F\|_{\mathcal N(D;H)} ,
\]
so that $\|F\ominus\tilde F\|_{L_\infty(D;G)}=\sup_x\dG(F(x),\tilde F(x))$ and $\|F\ominus\tilde F\|_{L_2(D;G)}^2=\int_D\dG(F(x),\tilde F(x))^2\,dx$.

The following result describes the \textit{transference principle}, where the approximation operator $\mathcal A$ is transferred to $G$.
\begin{proposition}
\label{prop:transference}
Let $\mathcal A:\mathcal N(D;H)\to\mathcal N(D;H)$ be any map and define
\[
  \mathcal A^G F := \exp\circ\,\mathcal A(\log\circ F), \qquad F\in\mathcal N(D;G).
\]
Then for every $F\in\mathcal N(D;G)$,
\[
  \|F\ominus \mathcal A^G F\|_{\mathcal N(D;G)}
  =
  \big\|\,\log\circ F - \mathcal A(\log\circ F)\,\big\|_{\mathcal N(D;H)} .
\]
In particular, $\mathcal A^G$ converges in $\mathcal N (D;G)$ at exactly the rate at which $\mathcal A$ approximates the lifted function $\log\circ F$ in $\mathcal N(D;H)$, and $\mathcal A^G$ inherits linearity, idempotence and norm bounds from $\mathcal A$ in the transported structure. 
\end{proposition}

\begin{proof}
By definition $\log\circ(\mathcal A^G F)=\mathcal A(\log\circ F)$, so the identity is the definition of $\|\cdot\ominus\cdot\|_{\mathcal N(D;G)}$. The inherited properties follow because $\exp$ is a linear isometry (Proposition~\ref{prop:hilbert}).
\end{proof}

To apply this identity to geometrically defined data, we relate chart regularity to ambient regularity and compare $d_G$ with natural distances. The following lemma identifies the pulled-back and extrinsic smoothness classes on compact ranges.

\begin{lemma} \label{lem:moser}
Let $\sigma>d/2$, let $G$ be an embedded submanifold of $\R^p$, for
instance the open subset $\Spd_n\subset\Sym_n$ or a matrix Lie group,
and let $F:D\to G$ have range in a compact set $K\subset G$. Then $F\in H^\sigma(D;\R^p)$ if and only if $\log\circ F\in H^\sigma(D;H)$, and
\[
  \|\log\circ F\|_{H^\sigma(D;H)}\le C_K\,\big(1+\|F\|_{H^\sigma(D;\R^p)}\big),\qquad
  \|F\|_{H^\sigma(D;\R^p)}\le C_K\,\big(1+\|\log\circ F\|_{H^\sigma(D;H)}\big),
\]
with a constant $C_K$ depending on $K$, $\sigma$, $D$ and the fixed chart and embedding. The same holds for the classes $C^k$.
\end{lemma}
The proof is given in Appendix~\ref{app:lemma_proof}. The constant $C_K$ depends on $F$ through $K$: if $B=\log\circ F\in H^\sigma(D;H)$, then $B$ is bounded, $F(\overline D)=\exp(B(\overline D))$ is compact in $G$, and $C_K$ depends on $F$ only through $\|B\|_{L_\infty(D;H)}$. This is harmless, since the error estimates below are stated in the pulled-back norm $\|F\|_{H^\sigma(D;G)}=\|B\|_{H^\sigma(D;H)}$ and never involve $C_K$. Conversely, $F\in H^\sigma(D;\R^p)$ alone does not give $B\in H^\sigma(D;H)$: the closure of the range must stay in $G$. For $\Spd_n$ this means eigenvalues bounded below by some $\lambda_->0$ and above by some $\lambda_+$, and the scalar case $F=e^B$, $F'=e^BB'$, shows how the two bounds enter: $\|B'\|_{L_2}\le\lambda_-^{-1}\|F'\|_{L_2}$ and $\|F'\|_{L_2}\le\lambda_+\|B'\|_{L_2}$. For example, $F(x)=x$ on $(0,1)$ is smooth with $\lambda_+=1$, but $\lambda_-=0$ and $\log F$ is unbounded.
\subsection{Positive definite matrices}

We write $\|\cdot\|_F$ for the Frobenius norm, $\|\cdot\|_2$ for the
spectral norm, and use the two increasing functions 
\begin{equation}
\label{eq:LamTheta}
  \Lam(\omega):=\frac{\sinh(\omega/2)}{\omega/2},\qquad
  \Th(a):=\frac{e^{a}-1}{a},\qquad \omega,a\ge0,
\end{equation}
with $\Lam(0)=\Th(0)=1$. For $A\in\Spd_n$,
\[
  \omega(A):=\lambda_{\max}(\log A)-\lambda_{\min}(\log A)=\log\cond_2(A)
\]
is the \emph{logarithmic condition number} of $A$; it vanishes exactly
for multiples of the identity and is invariant under scaling $A\mapsto
tA$. The affine-invariant metric $\ip{V}{W}_P=\tr(P^{-1}VP^{-1}W)$ on
$T_P\Spd_n=\Sym_n$ has the norm $\|V\|_P=\|P^{-1/2}VP^{-1/2}\|_F$ and
the geodesic distance 
\begin{equation}
\label{eq:dAI}
  \dAI(A,B)=\big\|\log\big(A^{-1/2}BA^{-1/2}\big)\big\|_F,
\end{equation}
see \cite[Ch.~6]{Bhatia2007}, \cite{Moakher2005,Pennec2006}.

\begin{theorem}[Three metrics on $\Spd_n$]
\label{thm:spd_metrics}
Assume $n\ge 2$ and let $A,B\in\Spd_n$, $X=\log A$, $Y=\log B$, and
$\omega(A,B):=\max\{\omega(A),\omega(B)\}$. 
\begin{enumerate}[label=\textup{(\roman*)},leftmargin=*]
  \item \textup{(Log-Euclidean versus affine-invariant.)}
  \[
    \|X-Y\|_F\;\le\;\dAI(A,B)\;\le\;\Lam(\omega(A,B))\,\|X-Y\|_F .
  \]
  Both inequalities are sharp: the left one is an equality when $A$ and $B$ commute, and
  \[
    \sup\Big\{\frac{\dAI(\tilde A,\tilde B)}{\|\log \tilde A-\log
        \tilde B\|_F}:\ \tilde A\ne
    \tilde B,\ \omega(\tilde A),\omega(\tilde B)\le\omega(A,B)\Big\}=\Lam(\omega(A,B)). 
  \]
  \item \textup{(Log-Euclidean versus Frobenius.)}
  \[
    \|A-B\|_F\;\le\;\max\{\|A\|_2,\|B\|_2\}\,\|X-Y\|_F,
    \qquad
    \|X-Y\|_F\;\le\;\frac{\|A-B\|_F}{\min\{\lambda_{\min}(A),\lambda_{\min}(B)\}} .
  \]
\end{enumerate}
\end{theorem}

The proof is given in Appendix~\ref{app:metrics}. The lower bound in (i) is the exponential metric increasing property \cite[Thm.~6.1.4]{Bhatia2007}, used in \cite{WendlandSPD}; it says that the log-Euclidean distance always underestimates the affine-invariant one, but does not by itself provide an upper estimate for the affine-invariant error. The upper bound supplies such an estimate, with $\Lam$ depending only on anisotropy, not scale. For $\omega\le2$, that is, condition numbers up to $e^2\approx7.4$, one has $\Lam\le1.18$, and for $\omega\le4$ (condition numbers up to $55$) $\Lam\le1.82$. That the two distances agree to first order at the identity and are equivalent on compact subsets is well known \cite{Arsigny2007,Pennec2006}; we are not aware of this sharp bound in terms of endpoint condition numbers in the literature.

\subsection{Lie groups}

Let $G$ be a Lie group with Lie algebra $\g$, let $\ip{\cdot}{\cdot}$ be an inner product on $\g$, and let $\dL$ be the distance of the corresponding left-invariant Riemannian metric; $(G,\dL)$ is complete. Set
\[
  c_\g:=\sup_{\|X\|=1}\|\ad_X\|_{\g\to\g},
\]
so that $\|[X,Y]\|\le c_\g\|X\|\|Y\|$, with $c_\g=0$ exactly for Abelian $\g$.

\begin{theorem}[Log-Euclidean versus invariant metrics]
\label{thm:lie_metrics}
Let $X,Y\in\g$.
\begin{enumerate}[label=\textup{(\roman*)},leftmargin=*]
  \item If $\|X\|,\|Y\|\le R$, then $\dL(\exp X,\exp Y)\le\Th(c_\g R)\,\|X-Y\|$.
  \item If the metric is bi-invariant (every compact Lie group admits such a metric), then 
  \[ \dL(\exp X,\exp Y)\le\|X-Y\| ,\quad   X,Y\in\g. \]
  \item If $G$ is exponential, then for every $R>0$ there is $c_R>0$ such that,
  \[
   \|X\|,\|Y\|\le R \quad  \Rightarrow \quad   c_R\,\|X-Y\|\le\dL(\exp X,\exp Y) .
  \]
\end{enumerate}
\end{theorem}

The proof is in Appendix~\ref{app:metrics}. Part (ii) does not use injectivity of $\exp$, which is what makes it applicable to compact groups in Section~\ref{sec:local}; for nilpotent groups the constant in (iii) can be made explicit, since $\ad_Z$ is nilpotent.

\subsection{Riemannian manifolds}

The comparison for Hadamard manifolds, and a local version for arbitrary manifolds, follow from Rauch's comparison theorem. Let $(M,g)$ be a complete Riemannian manifold with distance $d_M$, let $p\in M$, and let $0<\rho\le\inj(p)$ (with $\rho=\infty$ allowed if $\inj(p)=\infty$), so that $\exp_p$ is a diffeomorphism from the ball $B(0,\rho)\subset T_pM$ onto $B_M(p,\rho)$ with inverse $\log_p$. For $\kappa\ge0$ put $S_\kappa(r):=\sinh(\sqrt\kappa\,r)/(\sqrt\kappa\,r)$. All curvature quotients below take their continuous value $1$ at zero.

\begin{theorem}[Local and Hadamard comparison]
\label{thm:riemann}
Let the sectional curvatures on $B_M(p,\rho)$ satisfy $-\kappa_-\le K\le\kappa_+$ with $\kappa_\pm\ge0$, and let $v,w\in T_pM$ with $\|v\|,\|w\|\le r<\rho$.
\begin{enumerate}[label=\textup{(\roman*)},leftmargin=*]
  \item $d_M(\exp_pv,\exp_pw)\le S_{\kappa_-}(r)\,\|v-w\|$.
  \item If $\sqrt{\kappa_+}\,r<\pi$ and $B_M(p,r)$ is geodesically convex, then
  \[
    \|v-w\|\le\frac{\sqrt{\kappa_+}\,r}{\sin(\sqrt{\kappa_+}\,r)}\,d_M(\exp_pv,\exp_pw).
  \]
\end{enumerate}
In particular, if $M$ is a Hadamard manifold with $K\ge-\kappa$, then
\[
  \|v-w\|\le d_M(\exp_pv,\exp_pw)\le S_\kappa(r)\,\|v-w\|
  \qquad\text{for all } v,w \text{ with } \|v\|,\|w\|\le r,
\]
and every $r>0$.
\end{theorem}

For $\Spd_n$ with the affine-invariant metric the sectional curvatures lie in $[-\frac12,0]$, and the theorem gives the constant $S_{1/2}(r)=\sinh(r/\sqrt2)/(r/\sqrt2)$ for $\|X\|_F,\|Y\|_F\le r$; since $\omega(e^X)\le\sqrt2\|X\|_F$, Theorem~\ref{thm:spd_metrics} is at least as good and usually much better, its constant depending on the anisotropy of the data only. For a compact Lie group with a bi-invariant metric, $\kappa_-=0$ and $\kappa_+=c_\g^2/4$, so part (i) recovers Theorem~\ref{thm:lie_metrics}(ii) inside the injectivity radius, whereas the latter needs no radius restriction.

\subsection{Transfer of error estimates}

The following corollary is the form in which the comparison theorems are used.

\begin{corollary}[Error transfer to natural metrics]
\label{cor:transfer_natural}
Let $F,\tilde F:D\to G$, assume that $\log\circ F$ is bounded, and suppose that $\sup_{x\in D}\dG(F(x),\tilde F(x))\le\varepsilon<\infty$.
\begin{enumerate}[label=\textup{(\alph*)},leftmargin=*]
  \item If $G=\Spd_n$, $\omega_F:=\sup_x\omega(F(x))$ and $R_F:=\sup_x\|\log F(x)\|_2$, then
  \[
    \sup_x\dAI\big(F(x),\tilde F(x)\big)\le\Lam(\omega_F+2\varepsilon)\,\varepsilon,
    \qquad
    \sup_x\big\|F(x)-\tilde F(x)\big\|_F\le e^{R_F+\varepsilon}\,\varepsilon .
  \]
  \item If $G$ is an exponential Lie group with a left-invariant metric and $R_F:=\sup_x\|\log F(x)\|$, then
  \[
    \sup_x\dL(F(x),\tilde F(x))\le\Th(c_\g(R_F+\varepsilon))\,\varepsilon .
  \]
  \item If $G=M$ is a Hadamard manifold with $K\ge-\kappa$ and $R_F:=\sup_x\|\log_pF(x)\|$, then $\sup_xd_M(F(x),\tilde F(x))\le S_\kappa(R_F+\varepsilon)\,\varepsilon$.
\end{enumerate}

The same constants give, for $1\le q\le\infty$,
\[
 \|d_{\rm nat}(F,\tilde F)\|_{L_q(D)}
 \le C(\varepsilon)\,\|\log\circ F-\log\circ\tilde F\|_{L_q(D;H)},
\]
with the analogous statement for the ambient Frobenius error. Thus, rates transfer whenever the logarithmic errors have a common supremum bound. For compact groups with a bi-invariant metric, where only a principal logarithm is available, Theorem~\ref{thm:compact} gives the transfer with constant one and without this bound.
\end{corollary}

\begin{proof}
(a) By Weyl's inequality \cite[Ch.~III]{Bhatia1997}, the extreme eigenvalues of $\log\tilde F(x)=\log F(x)+E$ with $\|E\|_2\le\|E\|_F\le\varepsilon$ differ from those of $\log F(x)$ by at most $\varepsilon$, so $\omega(\tilde F(x))\le\omega_F+2\varepsilon$ and $\|\log\tilde F(x)\|_2\le R_F+\varepsilon$; apply Theorem~\ref{thm:spd_metrics} pointwise. (b) and (c) are Theorems~\ref{thm:lie_metrics}(i) and~\ref{thm:riemann}(i) with $R=R_F+\varepsilon$, respectively $r=R_F+\varepsilon$. The $L_q$ statement follows by applying the pointwise bounds with these uniform constants and integrating.
\end{proof}

\section{Single-scale kernel approximation}
\label{sec:singlescale}

We use kernels $K(x,y)=k(x,y)I_H$, where $k$ is a continuous scalar strictly positive definite kernel on $\overline D\times\overline D$. Its vector-valued native space $\Hk(D;H)$ has squared norm equal to the sum of squared scalar native norms in an orthonormal basis of $H$ \cite{MicchelliPontil,CarmeliDeVitoToigo2006,WendlandSPD}. For distinct sites $X=\{x_1,\ldots,x_N\}\subset\overline D$, write $V_X=\{\sum_i k(\cdot,x_i)\beta_i:\beta_i\in H\}$. Point values on the boundary use the continuous representative.

\begin{theorem}
\label{thm:kernel_LS}
For data $Y=(y_1,\ldots,y_N)\in H^N$, the interpolant $s_{0,Y}\in V_X$ exists uniquely and has minimal native norm among all interpolants. If $Y=u|_X$, it equals the native-space orthogonal projection $P_Xu$ onto $V_X$. For $\lambda>0$, the minimiser over $\Hk(D;H)$ of
\[
 \sum_i\|s(x_i)-y_i\|_H^2+\lambda\|s\|_{\Hk}^2
\]
is $s_{\lambda,Y}\in V_X$. In both cases the coefficients solve $(\mathbf K+\lambda I)\beta=Y$, where $\mathbf K=(k(x_i,x_j))_{i,j}$ and $\lambda=0$ for interpolation.
\end{theorem}

These are the componentwise representer and projection results of \cite[Thms.~4.1 and~4.6]{WendlandSPD}. One scalar matrix serves all $m$ components, and every approximant has the form
\begin{equation}\label{eq:scalar_form}
 \mathcal A(Y)(x)=\sum_i a_i(x)y_i,
\end{equation}
with coefficients independent of the data.

\begin{assumption}
\label{ass:kernel}
$\Phi\in L_1(\R^d)\cap C(\R^d)$ is even and its Fourier transform satisfies, for some $\sigma>d/2$ and $0<c_1\le c_2$,
\begin{equation}
\label{eq:fourier}
  c_1\,(1+\|\omega\|_2^2)^{-\sigma}\;\le\;\widehat\Phi(\omega)\;\le\;c_2\,(1+\|\omega\|_2^2)^{-\sigma},
  \qquad\omega\in\R^d.
\end{equation}
\end{assumption}

Under Assumption~\ref{ass:kernel}, $k(x,y):=\Phi(x-y)$ is positive
definite and its native space $\Nat\Phi(D)$ on the Lipschitz domain
$D$ coincides with $H^\sigma(D)$ with equivalent norms
\cite[Ch.~10]{Wendland2005}; we write $\Ceq$ for a constant with
$\Ceq^{-1}\|u\|_{H^\sigma(D)}\le\|u\|_{\Nat\Phi(D)}\le\Ceq\|u\|_{H^\sigma(D)}$. Examples
are the Mat\'ern kernels and the compactly supported functions
$\phi_{d,\ell}$ of \cite[Ch.~9]{Wendland2005} with $\ell\ge1$, for which
$\sigma=d/2+\ell+1/2$; the latter are used from
Section~\ref{sec:multiscale} on. The fill distance and the separation radius of a finite set $X\subset\overline D$ are
\[
  h_X:=\sup_{x\in D}\min_{x_j\in X}\|x-x_j\|_2,\qquad
  q_X:=\tfrac12\min_{i\ne j}\|x_i-x_j\|_2 .
\]
The error analysis rests on sampling inequalities. There are constants $C_S$ and $h_0>0$, depending only on $D$, $\sigma$, $d$ and $m$, such that for every finite $X\subset\overline D$ with $h_X\le h_0$ and every $u\in H^\sigma(D;H)$,
\begin{align}
\label{eq:SI_zero}
  &\|u\|_{L_2(D;H)}\le C_S\,h_X^{\sigma}\,\|u\|_{H^\sigma(D;H)},\qquad
  \|u\|_{L_\infty(D;H)}\le C_S\,h_X^{\sigma-d/2}\,\|u\|_{H^\sigma(D;H)}
  &&\text{if }u|_X=0,\\
\label{eq:SI_gen}
  &\|u\|_{L_\infty(D;H)}\le C_S\big(h_X^{\sigma-d/2}\,\|u\|_{H^\sigma(D;H)}+\max_i\|u(x_i)\|_H\big).
\end{align}
In the scalar case these are due to \cite{NWW2005,WendlandRieger2005,Arcangeli2007}, see also \cite[Ch.~11]{Wendland2005}; the $H$-valued versions, including corresponding $W^\tau_q$-bounds, are \cite[Thm.~4.3]{WendlandSPD}, stated there for $H=\Sym_n$ with a proof that applies verbatim to any finite-dimensional $H$.

Now let $G$ be a log-Hilbert target, $K=\Phi(\cdot-\cdot)\,I_H$, and let $F:D\to G$ be sampled at $X$. With $b_i:=\log F(x_i)$ and $Y=(b_1,\dots,b_N)$ we define
\[
  I^G_XF:=\exp\circ\,s_{0,Y},\qquad Q^G_{X,\lambda}F:=\exp\circ\,s_{\lambda,Y}.
\]
These are $G$-valued by construction and as smooth as the kernel.

\begin{theorem}[Single-scale errors] \label{thm:single_scale}
Let $F\in H^\sigma(D;G)$ and $h_X\le h_0$.
\begin{enumerate}[label=\textup{(\roman*)},leftmargin=*]
  \item \textup{(Interpolation.)}
  \begin{align*}
    \sup_{x\in D}\dG\big(F(x),I^G_XF(x)\big)&\le C\,h_X^{\,\sigma-d/2}\,\|F\|_{H^\sigma(D;G)},\\
    \|F\ominus I_X^GF\|_{L_2(D;G)}&\le C\,h_X^\sigma\|F\|_{H^\sigma(D;G)},
  \end{align*}
  \item \textup{(Penalized least squares.)} For $\lambda>0$,
  \[
    \sup_{x\in D}\dG\big(F(x),Q^G_{X,\lambda}F(x)\big)
    \le C\big(h_X^{\,\sigma-d/2}+\sqrt{\lambda}\,\big)\|F\|_{H^\sigma(D;G)} .
  \]
\end{enumerate}
The constants depend only on $C_S$ and $\Ceq$. 
\end{theorem}

\begin{proof}
Since $\log\circ I^G_XF$ and $\log\circ Q^G_{X,\lambda}F$ are the interpolant and the penalized least-squares approximant of $B:=\log\circ F$, Proposition~\ref{prop:transference} identifies the left-hand sides with the errors of the lifted problems, and these are bounded in \cite[Thms.~4.8 and~4.10]{WendlandSPD}. The proofs there use only the sampling inequalities, the minimal-norm property of Theorem~\ref{thm:kernel_LS} and the norm equivalence, so they hold for any finite-dimensional $H$.
\end{proof}

As for Theorem~\ref{thm:single_scale},  we note that Corollary~\ref{cor:transfer_natural} transfers the displayed errors to natural metrics when the logarithmic errors are uniformly bounded; for interpolation this follows from the displayed supremum estimate. In addition, for $\Spd_n$-valued data $A_i$, the computation solves $\sum_j(k(x_i,x_j)+\lambda\delta_{ij})C_j=\log A_i$ for $C_j\in\Sym_n$ and evaluates $\exp(\sum_jk(x,x_j)C_j)$. Therefore, there is one logarithm per distinct input and one exponential per output.

\paragraph{Preservation and equivariance.}
Any operator of the form \eqref{eq:scalar_form}, including the multilevel operator below, preserves linear subspaces of $H$. Unit-determinant positive definite data therefore give unit-determinant approximants, since $\operatorname{tr}(\log A_i)=0$. It also commutes with every log-linear map $\psi$ satisfying $\log\circ\psi=L\circ\log$, $L\in\mathcal L(H)$. Examples are orthogonal conjugation, inversion and powers on $\Spd_n$, and conjugation and inversion on exponential Lie groups.

For a log-affine map $\log\circ\psi=A+L\circ\log$, the defect is
\[
 \mathcal A(A+LB)-(A+L\mathcal AB)=(\mathcal A1-1)A.
\]
Hence constant reproduction gives affine equivariance. MLS reproduces constants; ordinary kernel interpolation need not. This explains the condition for respecting $P\mapsto tP$ on $\Spd_n$. Left translations of a general Lie group need not be log-affine, but in every two-step nilpotent group
\[
 \log(\exp A\exp X)=A+(I+\tfrac12\ad_A)X.
\]
Constant reproduction therefore gives left equivariance on these groups, including the Heisenberg group. General congruences of positive definite matrices need not be log-affine, so their equivariance is not automatic.

\section{Multilevel residual correction} \label{sec:multiscale}

We work in the model space and lift $F$ to $B=\log\circ F$. The multilevel scheme of \cite[Alg.~4.1]{sharon2022multiscale} then reads as follows.

\begin{algorithm}[H]
\caption{Multilevel approximation in $H$ (residual correction)}
\label{alg:multiscale_H}
\begin{enumerate}[leftmargin=1.8em,nosep,label=\arabic*.]
  \item Input: sampled $B$, sites $X_1,\ldots,X_L$, and linear operators $Q_j$
  using only samples on $X_j$.
  \item Set $B_0=0$ and $e_0=B$.
  \item For $j=1,\ldots,L$, set $s_j=Q_je_{j-1}$,
  $B_j=B_{j-1}+s_j$, and $e_j=e_{j-1}-s_j$.
  \item Return $B_L$ and $F_L=\exp\circ B_L$.
\end{enumerate}
\end{algorithm}

Since the $Q_j$ are linear,
\begin{equation}
\label{eq:recursion}
  e_L=(I-Q_L)(I-Q_{L-1})\cdots(I-Q_1)B,\qquad B_L=B-e_L ,
\end{equation}
and by Proposition~\ref{prop:transference}, $\|F\ominus F_L\|_{\mathcal N(D;G)}=\|e_L\|_{\mathcal N(D;H)}$ for every norm $\mathcal N$. In other words, $B_L=\mathcal M_LB$ with the linear multilevel operator
\begin{equation}
\label{eq:ML}
  \mathcal M_L:=I-(I-Q_L)(I-Q_{L-1})\cdots(I-Q_1),
\end{equation}
the rightmost factor being applied first.

Write $\Llog F=\log\circ F$, $\Eexp B=\exp\circ B$, and $Q_j^G=\Eexp\circ Q_j\circ\Llog$. In a function space $\mathcal N(D;H)$ preserved by every $Q_j$, the corresponding recursion on $G$ is
\begin{equation}\label{eq:recursion_G}
 F_0\equiv e,\quad R_{j-1}=F\ominus F_{j-1},\quad
 S_j=Q_j^GR_{j-1},\quad F_j=F_{j-1}\oplus S_j.
\end{equation}

\begin{proposition}[Multilevel conjugacy]\label{prop:conjugacy}
For $B=\Llog F$, the iterates of \eqref{eq:recursion_G} and Algorithm~\ref{alg:multiscale_H} satisfy
\[
 \Llog R_{j-1}=e_{j-1},\qquad \Llog S_j=s_j,\qquad \Llog F_j=B_j.
\]
Consequently $\mathcal M_L^G=\Eexp\circ\mathcal M_L\circ\Llog$, and
\[
\begin{tikzcd}[column sep=huge,row sep=large]
\mathcal N(D;G)\arrow[r,"\mathcal M_L^G"]\arrow[d,"\Llog"']&\mathcal N(D;G)\\
\mathcal N(D;H)\arrow[r,"\mathcal M_L"']&\mathcal N(D;H)\arrow[u,"\Eexp"']
\end{tikzcd}
\]
commutes at every level.
\end{proposition}
\begin{proof}
Since $e_j=B-B_j$ and $\Llog F_0=0$, induction gives
\[
 \Llog R_{j-1}=B-B_{j-1},\quad
 \Llog S_j=Q_je_{j-1},\quad
 \Llog F_j=B_{j-1}+s_j=B_j.
\]
The operator identity follows from \eqref{eq:ML}.
\end{proof}

This algebraic identity uses the transported operations and does not require Lemma~\ref{lem:moser}. It does not identify the method with an intrinsic scheme using Riemannian logarithms or Karcher means. For a principal logarithm on $\exp(U)$, the same identities require the logarithms of all residuals, corrections, and iterates to remain in $U$. The Hilbert-space scheme and its exponential output remain defined without this restriction.

\subsection{Fixed kernels and nested sites}
\label{subsec:telescope}

Let $K=k\,I_H$ be fixed and let $P_j$ be the interpolation operator on $X_j$ in $\Hk(D;H)$, that is, the $\Hk$-orthogonal projection onto $V_j=\spann\{K(\cdot,x)\alpha:\ x\in X_j,\ \alpha\in H\}$ (Theorem~\ref{thm:kernel_LS}). Nested sites $X_1\subset\dots\subset X_L$ give nested subspaces $V_1\subset\dots\subset V_L$.

The following elementary observation explains why residual correction with a fixed kernel adds nothing and motivates the scaled kernels of Section~\ref{subsec:scaled}.
\begin{proposition}[Fixed kernel, nested sites]
\label{thm:telescope}
Let $Q_j=P_j$ be exact interpolation in a fixed native space and let $X_1\subset\cdots\subset X_L$. Then the iterates of Algorithm~\ref{alg:multiscale_H} satisfy
\[
  B_j=P_j B,\qquad e_j=(I-P_j)B,\qquad j=1,\dots,L,
\]
and therefore $F_L=\exp(P_L\log F)=I^G_{X_L}F$: the multilevel interpolant is identical to single-scale interpolation on the finest level, and the error bounds of Theorem~\ref{thm:single_scale}(i) hold with $h_X=h_{X_L}$.
\end{proposition}

\begin{proof}
By induction, $B_1=P_1B$. If $B_{j-1}=P_{j-1}B\in V_{j-1}\subset V_j$, then $B_j=B_{j-1}+P_j(B-B_{j-1})\in V_j$, and for $x\in X_j$, $B_j(x)=B_{j-1}(x)+(B(x)-B_{j-1}(x))=B(x)$. Thus $B_j\in V_j$ interpolates $B$ on $X_j$, and by uniqueness of the interpolant in $V_j$, $B_j=P_jB$.
\end{proof}

Thus a fixed kernel and nested sites give no change in approximation power. For level-dependent kernels, the scalar multilevel interpolant has $N_L$ cardinal basis functions and is a projection onto their span \cite{GollwitzerKempfWendland2026}; the corresponding $H$-valued range has dimension $mN_L$. Multilevel effects require changing kernels, operators, regularisation, or sampling sets. 

\subsection{Compactly supported kernels scaled to the fill distance}
\label{subsec:scaled}

The classical remedy \cite{FloaterIske1996, Wendland2010,
  LeGiaSloanWendland2010, TownsendWendland2013} is to use on level $j$
a compactly supported kernel whose support radius is proportional to
the fill distance $h_j$. Then the native space norms change from
level to level, the linear systems are sparse and uniformly
well-conditioned on quasi-uniform sites, and convergence has to be proved by a different
argument. For scalar functions this is \cite[Thm.~1]{Wendland2010} in
$L_2$ and, for $W^\tau_q$-norms including $L_\infty$,
\cite[Thm.~2.5]{GollwitzerKempfWendland2026}. The following theorem provides the $H$-valued generalization of these two results with an explicit contraction factor. Its proof, deferred to Appendix~\ref{app:multiscale}, closely follows \cite{Wendland2010} but is included to derive the constants differently.

For $\delta>0$ let $\Phi_\delta:=\delta^{-d}\Phi(\cdot/\delta)$, so that $\widehat{\Phi_\delta}(\omega)=\widehat\Phi(\delta\omega)$, and let $\Nat{\Phi_\delta}(D;H)$ be the native space of the kernel $\Phi_\delta I_H$ on $D$.

\begin{theorem}[Multilevel convergence with scaled kernels]
\label{thm:multiscale_scaled}
Let $\Phi$ satisfy Assumption~\ref{ass:kernel} with $\supp\Phi\subset\overline{B(0,1)}$, and let $B\in H^\sigma(D;H)$. Let $X_1,X_2,\dots\subset\overline D$ be finite sets, not necessarily nested, whose fill distances satisfy $h_j=\mu^{j-1}h_1$ for some $0<\mu<1$ and $h_1\le h_0$. Fix $\nu>1$ with $\nu h_1\le1$, set $\delta_j:=\nu h_j$, and let $Q_j=P_j$ be interpolation on $X_j$ with the kernel $\Phi_{\delta_j}I_H$ in $\Nat{\Phi_{\delta_j}}(D;H)$. Let $B_L$ be the output of Algorithm~\ref{alg:multiscale_H}. Then there is a constant $C_2\ge1$, depending only on $D$, $d$, $\sigma$, $c_1$ and $c_2$, such that, with
\[
  \alpha:=C_2\big(\mu^\sigma+C_S\,\nu^{-\sigma}\big),
\]
the following hold for all $L\ge1$.
\begin{enumerate}[label=\textup{(\roman*)},leftmargin=*]
  \item \textup{(Error bounds.)}
  \begin{align}
  \label{eq:ms_L2}
    \|B-B_L\|_{L_2(D;H)}&\le 2C_SC_2\,\nu^{-\sigma}\,\alpha^{\,L-1}\,\|B\|_{H^\sigma(D;H)},\\
  \label{eq:ms_Linf}
    \|B-B_L\|_{L_\infty(D;H)}&\le 2C_SC_2\,\nu^{-\sigma}\,h_L^{-d/2}\,\alpha^{\,L-1}\,\|B\|_{H^\sigma(D;H)}.
  \end{align}
  \item \textup{(Rates.)} If $\alpha<1$, then $\sigma_{\mathrm{eff}}:=\log\alpha/\log\mu$ satisfies $0<\sigma_{\mathrm{eff}}<\sigma$, and in terms of the finest fill distance
  \begin{equation}
  \label{eq:ms_order}
  \begin{aligned}
    \|B-B_L\|_{L_2(D;H)}&\le 2C_SC_2\,\nu^{-\sigma}\Big(\frac{h_L}{h_1}\Big)^{\sigma_{\mathrm{eff}}}\|B\|_{H^\sigma(D;H)},\\
    \|B-B_L\|_{L_\infty(D;H)}&\le 2C_SC_2\,\nu^{-\sigma}\,h_1^{-d/2}\Big(\frac{h_L}{h_1}\Big)^{\sigma_{\mathrm{eff}}-d/2}\|B\|_{H^\sigma(D;H)} .
  \end{aligned}
  \end{equation}
  The $L_\infty$ bound decreases with $L$ if and only if $\alpha<\mu^{d/2}$, that is, $\sigma_{\mathrm{eff}}>d/2$.
  \item \textup{(Manifold-valued functions.)} If $B=\log\circ F$ with $F\in H^\sigma(D;G)$ and $F_L:=\exp\circ B_L$, then
  \[
    \|B-B_L\|_{L_2(D;H)}=\|F\ominus F_L\|_{L_2(D;G)},\qquad
    \|B-B_L\|_{L_\infty(D;H)}=\sup_{x\in D}\dG\big(F(x),F_L(x)\big).
  \]
  If $\alpha<\mu^{d/2}$, the supremum is bounded uniformly in $L$, and Corollary~\ref{cor:transfer_natural} transfers the $L_2$ rate in \eqref{eq:ms_order} to the natural distances and to the ambient Frobenius error.
\end{enumerate}
\end{theorem}

\begin{remark}[Variants of Theorem~\ref{thm:multiscale_scaled}]
\label{rem:multiscale_scaled}
Three notes for the last theorem:
\begin{enumerate}[label=\textup{(\roman*)},leftmargin=*]
\item Since $C_SC_2\nu^{-\sigma}\le\alpha$, the bounds \eqref{eq:ms_L2}--\eqref{eq:ms_Linf} imply the simpler
\[
  \|B-B_L\|_{L_2(D;H)}\le 2\,\alpha^{L}\,\|B\|_{H^\sigma(D;H)},\qquad
  \|B-B_L\|_{L_\infty(D;H)}\le 2\,h_L^{-d/2}\alpha^{L}\,\|B\|_{H^\sigma(D;H)},
\]
which is the form of \cite[Thm.~1]{Wendland2010}. The bounds \eqref{eq:ms_L2}--\eqref{eq:ms_Linf} are sharper by the factor $C_SC_2\nu^{-\sigma}/\alpha<1$ and make the contribution of the first level visible, see (iii).
\item The proof uses the relations between the $h_j$ and $\delta_j$ only through $\delta_j\le1$, an upper bound on $h_{j-1}/\delta_j$, and $h_L\le\mu^{L-1}h_1$. Hence, as in \cite{Wendland2010}, the hypotheses $h_j=\mu^{j-1}h_1$ and $\delta_j=\nu h_j$ can be relaxed to
\[
  c\mu h_{j-1}\le h_j\le\mu h_{j-1},\qquad \nu h_j\le\delta_j\le1,
\]
with some $c\in(0,1]$, and the theorem holds with $\alpha$ replaced by $C_2\big(\mu^\sigma+C_S(c\nu)^{-\sigma}\big)$. Quasi-uniform sites and an upper bound $\delta_j\lesssim h_j$ are needed only for the sparsity and conditioning statements of Section~\ref{sec:computation}.
\item For $L=1$, \eqref{eq:ms_L2} reads $\|B-B_1\|_{L_2(D;H)}\le2C_SC_2(h_1/\delta_1)^\sigma\|B\|_{H^\sigma(D;H)}$: this single-scale estimate decreases only through the ratio $h_1/\delta_1$. Since $\delta_1\le1$, this ratio is at least $h_1$, so the bound is never better than the fixed-kernel order $h_1^\sigma$, and for $\delta_1=\nu h_1$ with $\nu$ fixed it does not decrease with $h_1$ at all. The bound does not establish the stagnation observed in Section~\ref{sec:numerics}; residual correction adds the factor $\alpha^{L-1}$.
\end{enumerate}
\end{remark}

The proven exponent satisfies
\[
 \sigma_{\rm eff}=\sigma-
 \frac{\log\!\big(C_2[1+C_S(\mu\nu)^{-\sigma}]\big)}{\log(1/\mu)}.
\]
It does not exceed the fixed-kernel $L_2$ order $\sigma$. With $\nu=\gamma/\mu$, the contraction factor takes the form $C\mu^\sigma$ used in \cite{Wendland2010}. The constants do not certify contraction for every practical parameter choice: without $\alpha<1$ the theorem asserts nothing about convergence, and under $\alpha<1$ alone only the $L_2$ rate is asserted. The numerical examples below test the observed behavior without certifying these conditions. Rough targets and penalized residual correction are treated in the scalar theory \cite[Thms.~2 and~3]{Wendland2010} and transfer componentwise.

\subsection{A remark on quasi-interpolation}
\label{subsec:generic}

In general, the single-scale operators need not be interpolation operators. Suppose the $Q_j$ are linear, uniformly stable in $C^k(\overline D;H)$ and of approximation order $h_j^k$ in $C^0$, that is, $\|u-Q_ju\|_{C^k}\le C_{k,k}\|u\|_{C^k}$ and $\|u-Q_ju\|_{C^0}\le C_{0,k}h_j^k\|u\|_{C^k}$ with constants independent of $j$. Then \eqref{eq:recursion} gives
\begin{align*}
  \sup_{x\in D}\dG\big(F(x),F_L(x)\big)=\|B-B_L\|_{C^0(\overline D;H)}
  &\le C_{0,k}\,h_L^k\,C_{k,k}^{\,L-1}\,\|F\|_{C^k(D;G)}\\
  &= C_{0,k}\,h_1^k\,(C_{k,k}\mu^k)^{L-1}\,\|F\|_{C^k(D;G)}
\end{align*}
for $h_j=\mu^{j-1}h_1$, which is the stability-type convergence result of \cite{sharon2022multiscale}. Moving least-squares operators satisfy these assumptions and act componentwise on $H$-valued data, so this result, along with the improved rates for MLS on regular grids in \cite{DurstWendland2026}, transfers to log-Hilbert targets word for word via Proposition~\ref{prop:transference}. The fixed-kernel analysis above does not provide such a uniform $C^k$-stability bound. For the scaled kernels of Section~\ref{subsec:scaled}, the multilevel operator is bounded from $H^\sigma$ to $L_\infty$ \cite[Cor.~2.6]{GollwitzerKempfWendland2026}, but whether its Lebesgue constant is bounded independently of the number of levels is open \cite[Sec.~3.2]{GollwitzerKempfWendland2026}.

\section{Computational effort}
\label{sec:computation}

Because $K=k\,I_H$, all linear algebra is scalar and shared by the $m$ components; the only $G$-specific cost is one $\log$ per data point and one $\exp$ per evaluation point. We compare the single-scale scheme with a fixed kernel on the finest set $X_L$ to the multilevel scheme with scaled kernels. The sets $X_j$ are assumed quasi-uniform, $h_j/q_j\le\cqu$ for all $j$ with $q_j$ the separation distance, and $h_j=\mu^{j-1}h_1$, $\delta_j=\nu h_j$. Then $N_j:=|X_j|$ satisfies $c\,h_j^{-d}\le N_j\le C\,h_j^{-d}$ with constants depending on $D$ and $\cqu$, so that
\begin{equation}
\label{eq:Nj}
  N_j\le C_{\mathrm{qu}}\,\mu^{d(L-j)}N_L,\qquad \sum_{j=1}^LN_j\le\frac{C_{\mathrm{qu}}}{1-\mu^d}\,N_L,
  \qquad L=1+\frac{\log(h_1/h_L)}{\log(1/\mu)}=\mathcal O(\log N_L).
\end{equation}

\begin{proposition}[Sparsity and level-independent conditioning]
\label{prop:cond}
Let $\Phi$ satisfy Assumption~\ref{ass:kernel} with $\supp\Phi\subset\overline{B(0,1)}$, let $X\subset\overline D$ have separation distance $q_X$, and let $\delta>0$. Then the interpolation matrix $A_{\delta,X}:=(\Phi_\delta(x_i-x_k))_{i,k}$ has at most $(1+\delta/q_X)^d$ nonzero entries per row, and
\[
  \cond_2(A_{\delta,X})=\cond_2\big(A_{\Phi,X/\delta}\big)\le\kappa\big(q_X/\delta\big),
\]
where $\kappa:(0,\infty)\to(0,\infty)$ is a decreasing function depending only on $\Phi$ and $d$. For $X=X_j$ and $\delta=\delta_j=\nu h_j$ with quasi-uniform sites, every row has at most $(1+\nu\cqu)^d$ nonzeros and $\cond_2(A_{\delta_j,X_j})\le\kappa(1/(\nu\cqu))$ independently of $j$.
\end{proposition}

\begin{proof}
The balls $B(x_k,q_X)$ are pairwise disjoint, and those with $\|x_i-x_k\|\le\delta$ lie in $B(x_i,\delta+q_X)$, which bounds their number by $((\delta+q_X)/q_X)^d$. Next, $A_{\delta,X}=\delta^{-d}A_{\Phi,Y}$ with $Y:=X/\delta$, a set with separation distance $q_Y=q_X/\delta$, so the condition numbers agree. For $Y$ one has $\lambda_{\max}(A_{\Phi,Y})\le\|A_{\Phi,Y}\|_\infty\le\Phi(0)(1+1/q_Y)^d$ by the first part, while the lower bound for the smallest eigenvalue of a kernel matrix \cite[Ch.~12]{Wendland2005} gives $\lambda_{\min}(A_{\Phi,Y})\ge c(q_Y)$ with a positive increasing function $c$ depending only on $\Phi$ and $d$ (under \eqref{eq:fourier}, $c(q)\sim q^{2\sigma-d}$ as $q\to0$). Hence $\kappa(q):=\Phi(0)(1+1/q)^d/c(q)$ works, and for quasi-uniform sites $q_{X_j}/\delta_j\ge1/(\nu\cqu)$.
\end{proof}

For a fixed kernel on quasi-uniform sites, the usual eigenvalue estimates instead give $\operatorname{cond}_2 A=\mathcal O(h^{-2\sigma})=\mathcal O(N^{2\sigma/d})$ \cite{Wendland2005,Schaback1995}. These are worst-case bounds. We report the measured condition numbers with the numerical results.

\begin{proposition}[Complexity]\label{prop:complexity}
For quasi-uniform sites with $h_j=\mu^{j-1}h_1$ and $\delta_j=\nu h_j$, residual correction evaluated at $M$ points can be implemented with
\[
 \mathcal O\!\left(m(1+\nu\cqu)^d(N_L+M)\log N_L\right)
\]
arithmetic operations and $\mathcal O((1+\nu\cqu)^dN_L+mN_L)$ memory, using CG to a fixed relative tolerance and bucket-based neighbour lookup in fixed dimension. It requires one logarithm per distinct site in $\bigcup_jX_j$, hence $\mathcal O(N_L)$ logarithms, and $M$ exponentials. There are exactly $N_L$ logarithms when the sets are nested. Constants depend on $\mu$, $\nu$, $\cqu$, the domain and the tolerance.
\end{proposition}

\begin{proof}
Set $K=(1+\nu\cqu)^d$. Each correction has at most $K$ active centers at a point. Evaluating the previous corrections on $X_j$ therefore costs $\mathcal O(jmKN_j)$, and evaluating their sum at $M$ points costs $\mathcal O(LmKM)$. Bucket structures with cell size proportional to $\delta_j$ find the active centers in $\mathcal O(K)$ work per query. Matrix assembly costs $\mathcal O(KN_j)$. By Proposition~\ref{prop:cond}, CG takes $\mathcal O(\sqrt\kappa\log(1/\mathrm{tol}))$ iterations with a level-independent bound on $\kappa$, hence $\mathcal O(mKN_j)$ work per level. Sum these costs using \eqref{eq:Nj}. The same sum bounds the stored coefficients, matrices, and number of distinct input sites.
\end{proof}

A dense fixed-kernel factorization costs $\mathcal O(N_L^3+mN_L^2)$ work and $\mathcal O(N_L^2+mN_L)$ memory. Unpreconditioned CG has the standard worst-case iteration bound $\mathcal O(N_L^{\sigma/d}\log(1/\mathrm{tol}))$, with $\mathcal O(mN_L^2)$ work per iteration; actual counts depend on the spectrum and right-hand sides. These cost estimates concern a fixed tolerance. Theorem~\ref{thm:multiscale_scaled} assumes exact solves, so a complete analysis of solver errors in the residual recursion is separate.

\section{Compact groups and Karcher means}
\label{sec:local}

The global theory excludes compact groups such as $\mathrm{SO}(p)$, which have no global logarithm. Because of the structure of bi-invariant metrics (Theorem~\ref{thm:lie_metrics}(ii)), considerably more than a local statement is nonetheless available.

\subsection{Groups with a principal logarithm}

Let $G$ be a Lie group with Lie algebra $\g$ carrying an inner product, and let $U\subset\g$ be an open neighborhood of $0$ such that $\exp|_U$ is a diffeomorphism onto the open set $\exp(U)\subset G$; we call $\log:=(\exp|_U)^{-1}$ the principal logarithm. For $G=\mathrm{SO}(p)$ one may take $U=\{X\in\mathfrak{so}(p):\|X\|_2<\pi\}$, whose image is the set of rotations without eigenvalue $-1$ \cite[Ch.~11]{Higham2008}; other Lie groups admit exponential charts near the identity. Suppose $F:D\to G$ takes values in $\exp(U)$ and $B:=\log\circ F\in H^\sigma(D;\g)$, which by Lemma~\ref{lem:moser} holds as soon as $F$ is extrinsically of class $H^\sigma$ with range in a compact subset of $\exp(U)$. Run any of the Hilbert-space schemes of Sections~\ref{sec:singlescale} and~\ref{sec:multiscale} on $B$ to obtain $B_L$, and set $F_L:=\exp\circ B_L$. Note that $F_L$ is always defined and $G$-valued, whether or not $B_L$ stays inside $U$.

\begin{theorem}[Compact groups] \label{thm:compact}
Let $G$ carry a bi-invariant Riemannian metric with distance $d$. Then, in the setting above and for every scheme,
\[
  \sup_{x\in D}d\big(F(x),F_L(x)\big)\le\|B-B_L\|_{L_\infty(D;\g)},
\]
so that every $L_\infty$-estimate of Theorems~\ref{thm:single_scale}, \ref{thm:multiscale_scaled} and of Section~\ref{subsec:generic} holds for the Riemannian error with constant one and without any condition on the approximant. If, in addition, $\|B-B_L\|_{L_\infty}<\dist(B(\overline D),\g\setminus U)$, then $B_L(D)\subset U$, $B_L=\log\circ F_L$, and for every compact $U'\subset U$ there is $c_{U'}>0$ with $c_{U'}\|X-Y\|\le d(\exp X,\exp Y)$ for $X,Y\in U'$, so that the estimate is two-sided on compact subsets of the chart.
\end{theorem}

\begin{proof}
The first statement is Theorem~\ref{thm:lie_metrics}(ii) applied pointwise; it does not use injectivity of $\exp$. If $B_L(D)\subset U$, then $\log\circ F_L=B_L$ by definition of the principal logarithm. The lower bound follows from the compactness argument in the proof of Theorem~\ref{thm:lie_metrics}(iii), applied to the smooth map $\log$ on the open set $\exp(U)\supset\exp(U')$.
\end{proof}

The computation behind Theorem~\ref{thm:lie_metrics}(ii) also gives $\|e^X-e^Y\|_F\le\|X-Y\|_F$ for $X,Y\in\mathfrak{so}(p)$, so the Frobenius error of rotation-valued approximants is likewise bounded by the Lie-algebra error with constant one; for groups without a bi-invariant metric, such as $\mathrm{SE}(p)$, Theorem~\ref{thm:lie_metrics}(i) gives the one-sided bound with the constant $\Th(c_\g R)$.

\subsection{Log-Euclidean versus Karcher means}

The manifold-valued scheme of \cite{sharon2022multiscale} replaces the linear average $\sum_iw_i(x)b_i$ underlying a quasi-interpolant by the Karcher mean of the $g_i$ with weights $w_i(x)$. On a Lie group with a bi-invariant metric, the log-Euclidean mean $\exp(\sum_iw_i\log g_i)$ used here differs from the Karcher mean by $\mathcal O(r^3)$ near the identity and by $\mathcal O(r^2)$ near a general point, where $r$ measures the spread of logarithms.

\begin{proposition}
\label{prop:karcher}
Let $G$ carry a bi-invariant Riemannian metric with distance $d$. There are $r_0>0$ and $C>0$, depending only on $(G,d)$, such that for all $0<r\le r_0$, all $g_i=\exp(X_i)$ with $\|X_i\|\le r$ ($i=1,\dots,N$), and all weights $w_i\ge0$ with $\sum_iw_i=1$, the Karcher mean $m_K$, the unique minimiser of $\phi(g):=\frac12\sum_iw_i\,d(g,g_i)^2$ in $B(e,2r)$, and the log-Euclidean mean $\bar g:=\exp(\bar X)$, $\bar X:=\sum_iw_iX_i$, satisfy
\[
  d(m_K,\bar g)\le C\,r^3 .
\]
\end{proposition}

The proof is given in Appendix~\ref{app:karcher}. The same argument applies to $\Spd_n$ with the affine-invariant metric, whose Riemannian logarithm is $\log_P(Q)=P^{1/2}\log(P^{-1/2}QP^{-1/2})P^{1/2}$, because the symmetric Baker--Campbell--Hausdorff expansion $\log(e^{A}e^{B}e^{A})=2A+B+\mathcal O((\|A\|+\|B\|)^3)$ has no second-order term, as in \cite{Arsigny2007}; there the Karcher mean is unique globally, and both means have the same determinant. The third order is a property of logarithms near zero and does not persist for data clustered around an arbitrary point, where the two means differ at second order.

\begin{proposition}
\label{prop:karcher_general}
Let $(G,d)$ be as in Proposition~\ref{prop:karcher}, or $G=\Spd_n$ with the affine-invariant metric, and let $R>0$. Use $H=\g$ in the Lie-group case and $H=\Sym_n$ for positive definite matrices. There are $r_1>0$, $c_R\ge1$ and $C_R>0$ such that for all $0<r\le r_1$, all $X_0\in H$ with $\|X_0\|\le R$, all $g_i=\exp(X_0+Y_i)$ with $\|Y_i\|\le r$ and all weights $w_i\ge0$ with $\sum_iw_i=1$, the Karcher mean $m_K$, the unique minimiser of $\phi$ in $B(\exp X_0,2c_Rr)$, and the log-Euclidean mean $\bar g=\exp(X_0+\sum_iw_iY_i)$ satisfy
\[
  d(m_K,\bar g)\le C_R\,r^2 .
\]
The exponent is attained: in $\Spd_2$ with $X_0=\diag(1,-1)$, $Y_\pm=\pm tE$, $E=E_{12}+E_{21}$, and equal weights, $\bar g=e^{X_0}$, while the Karcher mean is the geometric mean $g_+\#g_-$ and
\[
  \dAI\big(g_+\#g_-,\,e^{X_0}\big)=\frac{\sinh2-2}{2\sqrt2}\,t^2+\mathcal O(t^4).
\]
\end{proposition}

For nonnegative local weights and a $C^1$ lifted target, the logarithmic spread is $r\lesssim h_j\|F\|_{C^1(D;G)}$. Proposition~\ref{prop:karcher_general} therefore gives an $\mathcal O(h_j^2)$ difference between the two single-level means, and an $\mathcal O(h_j^3)$ bound where $\|\log F(x)\|\lesssim h_j$. Higher-order MLS may use negative weights, which are outside these propositions. A comparison propagated through the multilevel residual recursion remains open.

\section{Numerical illustration}
\label{sec:numerics}

We examine whether multilevel residual correction can combine convergence with the sparse, well-conditioned systems obtained from scaled, compactly supported kernels. The experiment compares three methods on the same matrix-valued function and sampling grids: multilevel interpolation with scaled kernels, single-scale interpolation with the same scaled kernels, and single-scale interpolation with a fixed kernel radius. This comparison isolates the effect of residual correction and illustrates the tradeoff between convergence and conditioning.

Our test function is a smooth field of $3\times3$ positive definite matrices on $D=(0,1)^2$. We construct it as $F=\exp\circ B$, where
\[
 B(x)=\begin{pmatrix}
 \sin(2\pi x_1)
 & \tfrac12\cos(\pi x_1x_2)
 & \tfrac3{10}\sin(\pi(x_1+x_2))\\
 \cdot & \cos(2\pi x_2) & \tfrac25x_1x_2\\
 \cdot & \cdot & x_1-x_2
 \end{pmatrix}.
\]
All three methods approximate the symmetric log-field $B$ and return to the positive definite cone through the matrix exponential. They therefore preserve positive definiteness by construction.

We use the Wendland $C^2$ kernel ($\sigma=5/2$)
$\Phi(x)=(1-\|x\|_2)_+^4(4\|x\|_2+1)$ and nested grids
\[
 X_j=\{0,2^{-j},\ldots,1\}^2,\qquad j=2,\ldots,7.
\]
These have $N_j=(2^j+1)^2$ sites and fill distance $h_j=2^{-j}/\sqrt{2}$. The scaled kernels have support radii $\delta_j=4\cdot2^{-j}$, corresponding to $\mu=1/2$ and $\nu=4\sqrt2$. The fixed-radius method uses $\delta=1/2$ and is computed through $j=6$.

\begin{figure}[tbp]
    \centering
    \includegraphics[width=0.7\linewidth]{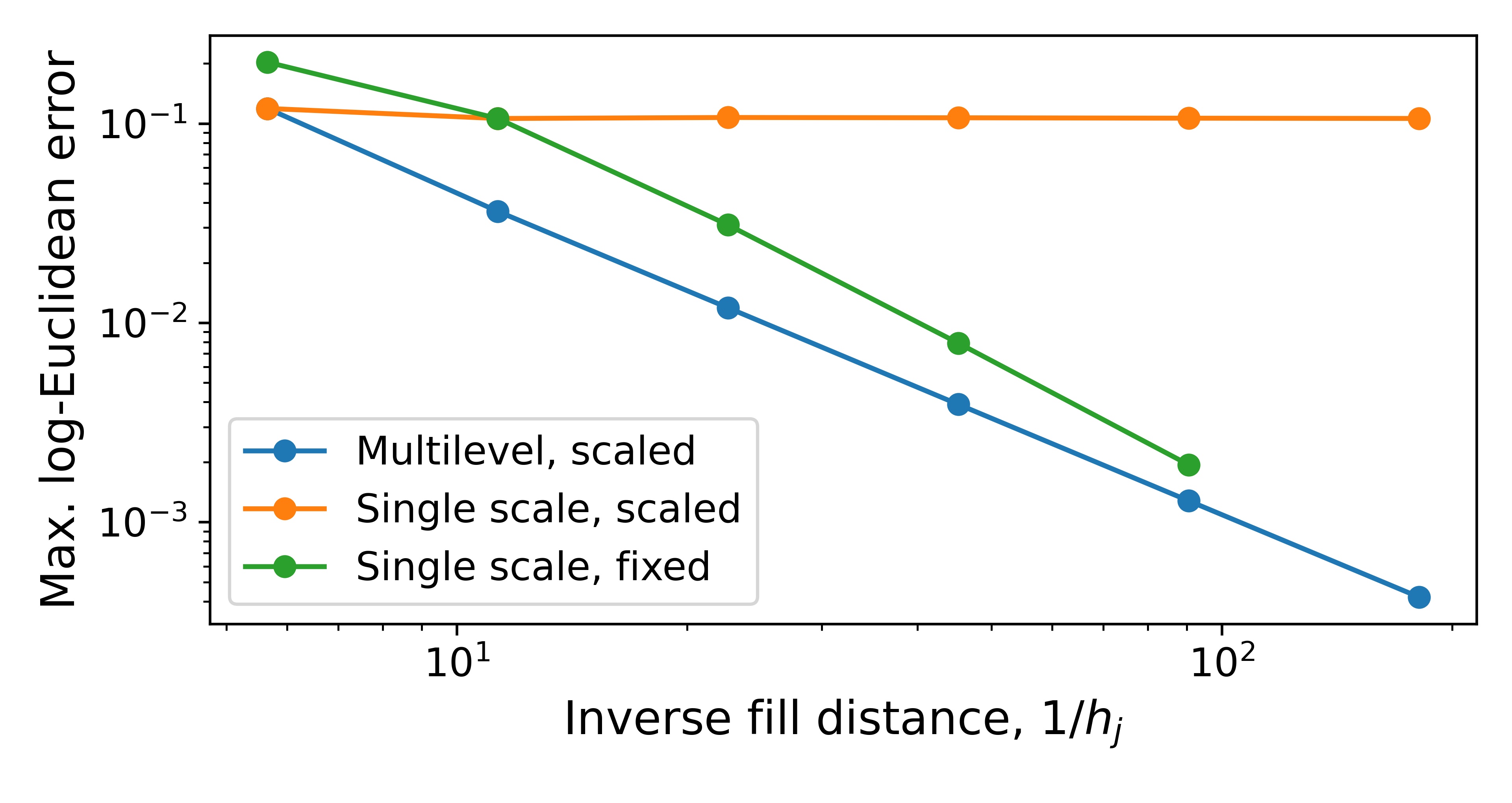}
    \caption{Comparison of multilevel interpolation with scaled kernels, single-scale interpolation with the same scaled kernels, and fixed-radius interpolation. The plotted errors are refined numerical estimates of the maximum log-Euclidean error. Residual correction produces decreasing errors while the corresponding single-scale scaled method stagnates. Fixed-radius interpolation also converges, but at the cost of solving systems with large condition numbers.}
    \label{fig:multiscale}
\end{figure}

Figure~\ref{fig:multiscale} shows the main effect. Single-scale interpolation with scaled kernels stagnates near an error of $0.1$ in this experiment, consistent with Remark~\ref{rem:multiscale_scaled}(iii), whereas applying those kernels to successive residuals steadily reduces the error. The multilevel error falls from $1.19\times10^{-1}$ to $4.21\times10^{-4}$, with a final measured dyadic rate of $1.60$. Thus, in this example, residual correction turns the stagnating single-scale construction into an effective approximation method.

Fixed-radius interpolation also converges, with errors comparable to those of the multilevel method. At the finest common grid, $N=4225$, the respective errors are $1.94\times10^{-3}$ and $1.28\times10^{-3}$. The distinction is in the linear systems: the fixed-radius condition number reaches $4.84\times10^6$, while the scaled systems remain below $195$ and have at most $45$ nonzeros per row. 

We validate the comparison in two ways. First, we check that the error estimates stabilize when evaluated at $256^2$ and $512^2$ cell midpoints, together with boundary points. Second, we check that the observed improvement also holds in the natural geometry. The affine-invariant and Frobenius errors, recorded in the accompanying data, exhibit comparable decay across the computed levels. Every evaluated pair satisfies the metric comparisons of Theorem~\ref{thm:spd_metrics} within floating-point tolerance. For the multilevel approximants, the largest sampled ratio $d_{\mathrm{AI}}/d_G$ is $1.169$, below the largest sampled distortion bound of approximately $1.229$. This supports the use of log-space approximation errors to describe accuracy in the natural metric.

As a separate implementation check, fixed-kernel multilevel interpolation and finest-level interpolation agree to $5.24\times10^{-15}$ in the sampled log-space maximum norm, as predicted by Proposition~\ref{thm:telescope}. The entire code and test appear in the accompanying public repository \url{https://github.com/nirsharon/log_Hilbert}.

\paragraph{Acknowledgements.}
NS is partially supported by the NSF-BSF award 2024791 and the BSF award 2024266. NS and HW are partially supported by the DFG award 514588180. 

\paragraph{Use of generative AI.}
ChatGPT (OpenAI) assisted with language editing and the
preparation of the Python numerical verification code and
a partial Lean~4 formalization. The authors retain full
responsibility for the manuscript, the code, and the
validation and interpretation of the reported results.

\paragraph{Code and data availability.}
The reproducibility and verification package is available at
\url{https://github.com/nirsharon/log_Hilbert}.
It contains the Python source, numerical data, and commands
for reproducing the experiments and the figure; no external
datasets are required. The Lean~4 development formalizes
the abstract results in Propositions~\ref{prop:hilbert},
\ref{prop:transference} and~\ref{prop:conjugacy}, together
with the nested orthogonal-projection argument underlying
Proposition~\ref{thm:telescope}. These formalizations were
mechanically checked by Lean~4. The analytic and geometric
error estimates are outside the scope of this formal
verification.

\bibliographystyle{plain}
\bibliography{refs}

\appendix

\section{Proof of Lemma~\ref{lem:moser}} \label{app:lemma_proof}

\begin{proof}
\emph{The composition estimate.}
For $\sigma>d/2$ the space $H^\sigma(D)$ is an algebra contained in $C(\overline D)$, and for a smooth map $\Psi$ between finite-dimensional Euclidean spaces the Moser-type composition estimate \cite[Ch.~13]{Taylor2011}, \cite[Ch.~5]{RunstSickel1996} states that $\Psi\circ u\in H^\sigma(D)$ for every $u\in H^\sigma(D)$, which is bounded since $\sigma>d/2$, with the tame bound
\[
  \|\Psi\circ u\|_{H^\sigma(D)}\le C\big(\Psi,\|u\|_{L_\infty(D)}\big)\,\big(1+\|u\|_{H^\sigma(D)}\big) .
\]
The point of the estimate is that its constant depends on $u$ only through $\|u\|_{L_\infty}$, and not on $\|u\|_{H^\sigma}$. It depends on $\Psi$ through finitely many of its derivatives on the closed ball of radius $\|u\|_{L_\infty}$, the precise number depending on the version invoked; the estimate of \cite[Ch.~13]{Taylor2011} uses a $C^N$ bound on $\Psi'$ with $N>\sigma$.

\emph{The first inequality.}
Let $\pi$ be the nearest-point projection onto $G$, which is smooth on a tubular neighbourhood of the compact set $K$, and let $\chi$ be a smooth cut-off function supported in that neighbourhood and equal to one near $K$. Then $\Psi:=\chi\cdot(\log\circ\pi)$ is smooth on $\R^p$, and $\Psi\circ F=\log\circ F$ because the range of $F$ lies in $K$. Apply the composition estimate to $u=F$: both $\|F\|_{L_\infty}$ and the derivatives of $\Psi$ on the relevant ball are bounded in terms of $K$, which gives the first inequality with a constant $C_K$.

\emph{The second inequality.}
Apply the composition estimate to $u=\log\circ F$, whose range lies in the compact set $\log(K)\subset H$, and to $\Psi=\chi'\cdot\exp$ with a cut-off $\chi'$ equal to one near $\log(K)$.

\emph{The classes $C^k$.}
The chain rule shows that $\Psi\circ u\in C^k(\overline D)$ whenever $u\in C^k(\overline D)$. For $k\ge1$, the linear bound $\|\Psi\circ u\|_{C^k}\le C(\Psi,\|u\|_{C^0})(1+\|u\|_{C^k})$ does not follow from the chain rule alone, which produces products of lower-order derivatives of $u$, but from the Gagliardo--Nirenberg interpolation inequalities
\[
  \|D^iu\|_{L_\infty}\le C\,\|u\|_{L_\infty}^{1-i/k}\,\|u\|_{C^k}^{i/k},\qquad 0\le i\le k,
\]
which bound every such product of total order $k$ by $C(\|u\|_{C^0})\,\|u\|_{C^k}$.
\end{proof}

\section{Proofs of the comparison theorems}
\label{app:metrics}

\subsection{Proof of Theorem~\ref{thm:spd_metrics}}

Both parts of the theorem are proved in the same way: we bound the differential of $\exp$, or of $\log$, along a curve and integrate. The differential is computed with the Dalecki\u{\i}--Kre\u{\i}n formula \cite[Thm.~V.3.3]{Bhatia1997}: if $f$ is $C^1$ on an open interval containing the spectrum of $Z=U\diag(z_1,\dots,z_n)U^{T}\in\Sym_n$, then
\begin{equation}
\label{eq:DK}
  \mathrm d f_Z[W]=U\big(f^{[1]}(z_i,z_j)\,\widetilde W_{ij}\big)_{i,j}U^{T},
  \quad \widetilde W:=U^{T}WU,\quad
  f^{[1]}(a,b):=\begin{cases}\dfrac{f(a)-f(b)}{a-b},& a\ne b,\\[1ex] f'(a),&a=b.\end{cases}
\end{equation}
Note that $\|\widetilde W\|_F=\|W\|_F$, since the Frobenius norm is unitarily invariant.

\begin{lemma}[Differentials of $\exp$ and $\log$]
\label{lem:dexp}
Let $Z\in\Sym_n$, $P\in\Spd_n$ and $W\in\Sym_n$. Then
\begin{align}
\label{eq:dexp_AI}
  \|W\|_F&\;\le\;\big\|\mathrm d\exp_Z[W]\big\|_{e^Z}\;\le\;\Lam\big(\omega(e^Z)\big)\,\|W\|_F,\\
\label{eq:dexp_F}
  \big\|\mathrm d\exp_Z[W]\big\|_F&\;\le\;e^{\lambda_{\max}(Z)}\,\|W\|_F,\\
\label{eq:dlog_F}
  \big\|\mathrm d\log_P[W]\big\|_F&\;\le\;\frac{\|W\|_F}{\lambda_{\min}(P)} .
\end{align}
\end{lemma}

\begin{proof}
For \eqref{eq:dexp_AI}, apply \eqref{eq:DK} to $f=\exp$ at $Z$ and conjugate by $e^{-Z/2}=U\diag(e^{-z_i/2})U^T$:
\[
  \big(U^T e^{-Z/2}\,\mathrm d\exp_Z[W]\,e^{-Z/2}U\big)_{ij}
  = e^{-(z_i+z_j)/2}\,\frac{e^{z_i}-e^{z_j}}{z_i-z_j}\,\widetilde W_{ij}
  = \Lam(z_i-z_j)\,\widetilde W_{ij},
\]
with $\Lam$ from \eqref{eq:LamTheta} extended evenly to negative arguments. Since $\|V\|_{e^Z}=\|e^{-Z/2}Ve^{-Z/2}\|_F$,
\[
  \big\|\mathrm d\exp_Z[W]\big\|_{e^Z}^2=\sum_{i,j}\Lam(z_i-z_j)^2\,\widetilde W_{ij}^2 ,
\]
and \eqref{eq:dexp_AI} follows from $1\le\Lam(z_i-z_j)\le\Lam(\lambda_{\max}(Z)-\lambda_{\min}(Z))=\Lam(\omega(e^Z))$, which holds because $\Lam$ is even and increasing on $[0,\infty)$.

For \eqref{eq:dexp_F}, the mean value theorem gives $|\exp^{[1]}(z_i,z_j)|\le e^{\max\{z_i,z_j\}}\le e^{\lambda_{\max}(Z)}$, so \eqref{eq:DK} yields $\|\mathrm d\exp_Z[W]\|_F\le e^{\lambda_{\max}(Z)}\|\widetilde W\|_F$.

For \eqref{eq:dlog_F}, apply \eqref{eq:DK} to $f=\log$ at $P$: for $a,b>0$ the mean value theorem gives $\log^{[1]}(a,b)\le1/\min\{a,b\}$, hence $\|\mathrm d\log_P[W]\|_F\le\|\widetilde W\|_F/\lambda_{\min}(P)$.
\end{proof}

\begin{proof}[Proof of Theorem~\ref{thm:spd_metrics}]
Throughout, $X=\log A$ and $Y=\log B$, and we abbreviate $\omega:=\omega(A,B)$.

\emph{Step 1: the upper bound in (i).}
Let $Z(t):=(1-t)X+tY$ and $\gamma(t):=e^{Z(t)}$, a smooth curve in $\Spd_n$ from $A$ to $B$ with $\gamma'(t)=\mathrm d\exp_{Z(t)}[Y-X]$. Since $\lambda_{\max}$ is convex and $\lambda_{\min}$ is concave on $\Sym_n$,
\[
  \omega\big(e^{Z(t)}\big)=\lambda_{\max}(Z(t))-\lambda_{\min}(Z(t))\le(1-t)\,\omega(A)+t\,\omega(B)\le\omega .
\]
Hence, by the upper bound in \eqref{eq:dexp_AI} and the monotonicity of $\Lam$,
\[
  \dAI(A,B)\le\int_0^1\big\|\mathrm d\exp_{Z(t)}[Y-X]\big\|_{\gamma(t)}\,dt
  \le \Lam(\omega)\,\|Y-X\|_F .
\]

\emph{Step 2: the lower bound in (i).}
Let $\gamma:[0,1]\to\Spd_n$ be any piecewise $C^1$ curve from $A$ to $B$, and let $\zeta:=\log\circ\gamma$, which is piecewise $C^1$ because $\log$ is smooth on $\Spd_n$. Then $\gamma'=\mathrm d\exp_\zeta[\zeta']$, and the lower bound in \eqref{eq:dexp_AI} gives $\|\gamma'(t)\|_{\gamma(t)}\ge\|\zeta'(t)\|_F$. Therefore the affine-invariant length of $\gamma$ satisfies
\[
  \int_0^1\|\gamma'(t)\|_{\gamma(t)}\,dt\;\ge\;\int_0^1\|\zeta'(t)\|_F\,dt\;\ge\;\|\zeta(1)-\zeta(0)\|_F=\|Y-X\|_F ,
\]
and taking the infimum over $\gamma$ gives $\dAI(A,B)\ge\|X-Y\|_F$. If $A$ and $B$ commute, they are simultaneously diagonalisable, $A^{-1/2}BA^{-1/2}=e^{Y-X}$, and \eqref{eq:dAI} reduces to $\|X-Y\|_F$.

\emph{Step 3: sharpness of $\Lam(\omega)$.}
It suffices to treat $n=2$; the general case follows by embedding as a leading block. For $\omega=0$, scalar multiples of the identity give equality. For $\omega>0$ and $0<\varepsilon<\omega/2$ let $W:=E_{12}+E_{21}$ and
\[
  a_\varepsilon:=\sqrt{\omega^2/4-\varepsilon^2},\qquad
  X_\varepsilon:=\diag(a_\varepsilon,-a_\varepsilon),\qquad
  Y_\varepsilon:=X_\varepsilon+\varepsilon W .
\]
Then $\omega(e^{X_\varepsilon})=2a_\varepsilon<\omega$, $\omega(e^{Y_\varepsilon})=2\sqrt{a_\varepsilon^2+\varepsilon^2}=\omega$, and $\|X_\varepsilon-Y_\varepsilon\|_F=\sqrt2\,\varepsilon$, so the pair $(e^{X_\varepsilon},e^{Y_\varepsilon})$ is admissible in the supremum. Along the smooth curve $s\mapsto e^{X_\varepsilon+sW}$, the first-order expansion of the Riemannian distance gives
\[
  \dAI\big(e^{X_\varepsilon},e^{Y_\varepsilon}\big)=\varepsilon\,\big\|\mathrm d\exp_{X_\varepsilon}[W]\big\|_{e^{X_\varepsilon}}+o(\varepsilon),
\]
and since $X_\varepsilon$ is diagonal, the computation in the proof of Lemma~\ref{lem:dexp} gives $\|\mathrm d\exp_{X_\varepsilon}[W]\|_{e^{X_\varepsilon}}=\sqrt2\,\Lam(2a_\varepsilon)$. Hence the ratio $\dAI(e^{X_\varepsilon},e^{Y_\varepsilon})/\|X_\varepsilon-Y_\varepsilon\|_F$ tends to $\Lam(\omega)$ as $\varepsilon\to0$, so the supremum is at least $\Lam(\omega)$; by Step~1 it is at most $\Lam(\omega)$.

\emph{Step 4: part (ii).}
Along $Z(t)$ from Step~1, convexity of $\lambda_{\max}$ gives
\[
  \lambda_{\max}(Z(t))\le\max\{\lambda_{\max}(X),\lambda_{\max}(Y)\},
  \qquad e^{\lambda_{\max}(X)}=\|A\|_2,\quad e^{\lambda_{\max}(Y)}=\|B\|_2 .
\]
Integrating \eqref{eq:dexp_F} along $\gamma$ therefore gives
\[
  \|A-B\|_F\le\int_0^1\big\|\mathrm d\exp_{Z(t)}[Y-X]\big\|_F\,dt\le\max\{\|A\|_2,\|B\|_2\}\,\|X-Y\|_F .
\]
For the second inequality, use the segment $P(t):=(1-t)A+tB$, which lies in $\Spd_n$ by convexity of the cone, with $\lambda_{\min}(P(t))\ge\min\{\lambda_{\min}(A),\lambda_{\min}(B)\}$ by concavity of $\lambda_{\min}$. Integrating \eqref{eq:dlog_F} along $P$ gives
\[
  \|X-Y\|_F\le\int_0^1\big\|\mathrm d\log_{P(t)}[B-A]\big\|_F\,dt\le\frac{\|A-B\|_F}{\min\{\lambda_{\min}(A),\lambda_{\min}(B)\}} .\qedhere
\]
\end{proof}

\subsection{Proof of Theorem~\ref{thm:lie_metrics}}

We use the formula for the differential of the exponential map \cite[Ch.~5]{Hall2015},
\begin{equation}
\label{eq:dexp_lie}
  \mathrm d\exp_X[W]
  = \mathrm dL_{\exp X}\Big(\frac{1-e^{-\ad_X}}{\ad_X}\,W\Big),
  \qquad
  \frac{1-e^{-\ad_X}}{\ad_X}:=\sum_{k\ge0}\frac{(-\ad_X)^k}{(k+1)!},
\end{equation}
where $\mathrm dL_g$, the differential of left translation, is an isometry from $\g=T_eG$ onto $T_gG$ for a left-invariant metric. As in the proof of Theorem~\ref{thm:spd_metrics}, the idea is to bound the differential along the segment $Z(t)=(1-t)X+tY$ and to integrate.

\begin{proof}[Proof of Theorem~\ref{thm:lie_metrics}]
\emph{Part (i).}
By \eqref{eq:dexp_lie} and the isometry property of $\mathrm dL_{\exp Z}$,
\[
  \|\mathrm d\exp_Z[W]\|_{\exp Z}\le\sum_{k\ge0}\frac{\|\ad_Z\|^k}{(k+1)!}\,\|W\|
  =\Th(\|\ad_Z\|)\,\|W\|\le\Th(c_\g\|Z\|)\,\|W\| ,
\]
where the last step uses $\|\ad_Z\|\le c_\g\|Z\|$ and the monotonicity of $\Th$. Along the segment $Z(t)=(1-t)X+tY$ one has $\|Z(t)\|\le R$, so the curve $t\mapsto\exp(Z(t))$ from $\exp X$ to $\exp Y$ has length at most $\Th(c_\g R)\|X-Y\|$.

\emph{Part (ii).}
For a bi-invariant metric the inner product on $\g$ is $\Ad$-invariant, so each $\ad_Z$ is skew-adjoint, hence normal with purely imaginary spectrum $\{i\theta_k\}$. Consequently $(1-e^{-\ad_Z})/\ad_Z$ is normal with eigenvalues
\[
  \frac{1-e^{-i\theta_k}}{i\theta_k}\quad\text{(equal to $1$ when $\theta_k=0$)},
  \qquad\text{of modulus}\qquad
  \frac{2\,|\sin(\theta_k/2)|}{|\theta_k|}\le1 .
\]
A normal operator has norm equal to its spectral radius, so $\|\mathrm d\exp_Z[W]\|_{\exp Z}\le\|W\|$ for all $Z\in\g$, and the claim follows by integrating along the segment as in (i). For $\mathrm{SO}(p)$ with the Frobenius inner product, left translation is also a Frobenius isometry of $\R^{p\times p}$, which gives the bound quoted in Section~\ref{sec:local}.

\emph{Part (iii).}
Let $K_R:=\exp(\overline{B}_R)$, a compact subset of $G$, and let $K_R':=\{g\in G:\dL(g,K_R)\le1\}$, which is compact by completeness of $(G,\dL)$ and the Hopf--Rinow theorem. Since $\log=\exp^{-1}$ is smooth on $G$,
\[
  \Gamma_R:=\sup_{g\in K'_R}\|\mathrm d\log_g\|<\infty .
\]
Let $g,h\in K_R$. If $\dL(g,h)<1$, a minimising geodesic from $g$ to $h$ exists and stays in $K'_R$, hence $\|\log g-\log h\|\le\Gamma_R\,\dL(g,h)$. If $\dL(g,h)\ge1$, then trivially $\|\log g-\log h\|\le2R\le2R\,\dL(g,h)$. Thus $c_R:=1/\max\{\Gamma_R,2R\}$ works.
\end{proof}

\subsection{Proof of Theorem~\ref{thm:riemann}}

\begin{proof}
\emph{Bounds on the differential of $\exp_p$.}
By Rauch's comparison theorem \cite[Ch.~10]{doCarmo1992}, for $u\in T_pM$ with $\|u\|=s<\rho$ and $\xi\perp u$,
\begin{equation}
\label{eq:rauch}
  \frac{\sin(\sqrt{\kappa_+}\,s)}{\sqrt{\kappa_+}\,s}\,\|\xi\|\;\le\;\|\mathrm d(\exp_p)_u[\xi]\|\;\le\;S_{\kappa_-}(s)\,\|\xi\| .
\end{equation}
The lower bound requires $\sqrt{\kappa_+}\,s<\pi$; the upper bound holds because there are no conjugate points along the geodesic $t\mapsto\exp_p(tu)$, $t\in[0,1]$, inside the injectivity radius. By the Gauss lemma, $\mathrm d(\exp_p)_u$ maps $u$ to a vector of norm $\|u\|$ orthogonal to the image of $u^\perp$, and since $\sin x/x\le1\le S_{\kappa_-}(s)$, both bounds in \eqref{eq:rauch} hold for all $\xi\in T_pM$.

\emph{Part (i).}
The segment $u(t):=v+t(w-v)$ stays in $\overline{B}(0,r)$, so $\exp_p\circ u$ is a curve from $\exp_pv$ to $\exp_pw$. By the upper bound in \eqref{eq:rauch} and the monotonicity of $S_{\kappa_-}$, its length is at most $S_{\kappa_-}(r)\|v-w\|$.

\emph{Part (ii).}
For $\|v\|,\|w\|<r$, geodesic convexity gives a minimising geodesic $\gamma$ from $\exp_pv$ to $\exp_pw$ in $B_M(p,r)$. Then $\zeta:=\log_p\circ\gamma$ is smooth, and \eqref{eq:rauch} gives
\[
  \|\zeta'\|\le\frac{\sqrt{\kappa_+}\,r}{\sin(\sqrt{\kappa_+}\,r)}\,\|\gamma'\| ,
\]
because $x\mapsto x/\sin x$ is increasing on $[0,\pi)$. Integrating gives the claim; the boundary cases follow by continuity.

\emph{The Hadamard case.}
On a Hadamard manifold $\inj(p)=\infty$, every ball is geodesically convex, and $\kappa_+=0$ gives the constant $1$ in (ii); hence both bounds hold for every $r>0$.
\end{proof}

\section{Proof of Theorem~\ref{thm:multiscale_scaled}}
\label{app:multiscale}

We use the Fourier transform $\widehat f(\omega)=\int_{\R^d}f(x)e^{-ix\cdot\omega}\,dx$ and write $a_d:=(2\pi)^{-d}$, so that the norms on $\R^d$ are
\[
 \|f\|_{H^\sigma(\R^d)}^2=a_d\int|\widehat f(\omega)|^2(1+\|\omega\|^2)^\sigma\,d\omega,
 \qquad
 \|f\|_{\Nat{\Phi_\delta}(\R^d)}^2=a_d\int\frac{|\widehat f(\omega)|^2}{\widehat\Phi(\delta\omega)}\,d\omega .
\]
On $D$ we use the restriction norms, $\|f\|_{\mathcal N(D)}:=\inf\{\|g\|_{\mathcal N(\R^d)}:g|_D=f\}$; the restricted native space is the reproducing kernel Hilbert space of the restricted kernel \cite[Ch.~10]{Wendland2005}. With these conventions, $\|f\|_{L_2(D)}\le\|f\|_{H^\sigma(D)}$.

The proof has two ingredients: the comparison between the native space norm of a scaled kernel and the Sobolev norm, which is the content of the following lemma, and the sampling inequalities \eqref{eq:SI_zero}.

\begin{lemma}[Native spaces of scaled kernels]
\label{lem:scaled_norms}
Let $\Phi$ satisfy Assumption~\ref{ass:kernel}. There is a constant $C_D$, depending only on $D$, $d$, $\sigma$ and $c_1$, such that for all $0<\delta\le1$ the sets $\Nat{\Phi_\delta}(D)$ and $H^\sigma(D)$ coincide and, for all scalar or $H$-valued $f$ with the same constants,
\begin{align}
\label{eq:scaled_a}
  \|f\|_{H^\sigma(D)}&\le\sqrt{c_2}\;\delta^{-\sigma}\,\|f\|_{\Nat{\Phi_\delta}(D)},\\
\label{eq:scaled_b}
  \|f\|_{\Nat{\Phi_\delta}(D)}&\le\sqrt{C_D}\,\big(\|f\|_{L_2(D)}+\delta^{\sigma}\|f\|_{H^\sigma(D)}\big).
\end{align}
\end{lemma}

\begin{proof}
\emph{Step 1: the native space norm on $\R^d$.}
Since $\widehat{\Phi_\delta}(\omega)=\widehat\Phi(\delta\omega)$, the bounds \eqref{eq:fourier} give, for every $f\in\Nat{\Phi_\delta}(\R^d)$,
\begin{equation}
\label{eq:native_fourier}
  \frac{1}{c_2}\int|\widehat f(\omega)|^2(1+\delta^2\|\omega\|^2)^\sigma\,d\omega
  \;\le\;(2\pi)^{d}\,\|f\|^2_{\Nat{\Phi_\delta}(\R^d)}
  \;\le\;\frac{1}{c_1}\int|\widehat f(\omega)|^2(1+\delta^2\|\omega\|^2)^\sigma\,d\omega .
\end{equation}

\emph{Step 2: proof of \eqref{eq:scaled_a}.}
For $\delta\le1$ one has $1+\delta^2\|\omega\|^2\ge\delta^2(1+\|\omega\|^2)$, so the lower bound in \eqref{eq:native_fourier} gives $\|f\|^2_{H^\sigma(\R^d)}\le c_2\,\delta^{-2\sigma}\|f\|^2_{\Nat{\Phi_\delta}(\R^d)}$, which is \eqref{eq:scaled_a} on $\R^d$ \cite[Lemma~1]{Wendland2010}. For $f\in\Nat{\Phi_\delta}(D)$ let $g$ be its minimal-norm extension to $\R^d$, so that $\|g\|_{\Nat{\Phi_\delta}(\R^d)}=\|f\|_{\Nat{\Phi_\delta}(D)}$. Then
\[
  \|f\|_{H^\sigma(D)}\le\|g\|_{H^\sigma(\R^d)}\le\sqrt{c_2}\,\delta^{-\sigma}\|g\|_{\Nat{\Phi_\delta}(\R^d)}=\sqrt{c_2}\,\delta^{-\sigma}\|f\|_{\Nat{\Phi_\delta}(D)} .
\]

\emph{Step 3: proof of \eqref{eq:scaled_b}.}
With $C_\sigma:=\max\{1,2^{\sigma-1}\}$ one has $(1+\delta^2\|\omega\|^2)^\sigma\le C_\sigma(1+\delta^{2\sigma}\|\omega\|^{2\sigma})$, so the upper bound in \eqref{eq:native_fourier} gives, for every $g\in H^\sigma(\R^d)$,
\[
  \|g\|^2_{\Nat{\Phi_\delta}(\R^d)}\le\frac{C_\sigma}{c_1}\Big(\|g\|^2_{L_2(\R^d)}+\delta^{2\sigma}\|g\|^2_{H^\sigma(\R^d)}\Big).
\]
Apply this to $g=Ef$, where $E$ is Stein's extension operator \cite[Ch.~VI]{Stein1970}, which for Lipschitz $D$ is bounded both from $L_2(D)$ to $L_2(\R^d)$ and from $H^\sigma(D)$ to $H^\sigma(\R^d)$. Since $\|f\|_{\Nat{\Phi_\delta}(D)}\le\|Ef\|_{\Nat{\Phi_\delta}(\R^d)}$, this gives \eqref{eq:scaled_b} with
\[
  C_D:=\frac{C_\sigma}{c_1}\max\Big\{\|E\|^2_{L_2(D)\to L_2(\R^d)},\ \|E\|^2_{H^\sigma(D)\to H^\sigma(\R^d)}\Big\}.
\]

\emph{Step 4: conclusion.}
Equality of the sets $\Nat{\Phi_\delta}(D)=H^\sigma(D)$ follows from the two inequalities, and the $H$-valued statements follow by Parseval's identity in an orthonormal basis of $H$.
\end{proof}

\begin{proof}[Proof of Theorem~\ref{thm:multiscale_scaled}]
Throughout, all norms are taken on $D$ and are $H$-valued. Write $\Nat j:=\Nat{\Phi_{\delta_j}}(D;H)$ and
\[
  C_2:=\max\big\{1,\sqrt{c_2C_D}\big\},\qquad \beta:=C_2\big(1+C_S(\mu\nu)^{-\sigma}\big),
\]
so that $\alpha=\mu^\sigma\beta$. Since $\delta_j\le\delta_1=\nu h_1\le1$ for all $j$, Lemma~\ref{lem:scaled_norms} applies on every level.

\emph{Step 1: one level of residual correction.}
Let $j\ge1$. Since $e_{j-1}\in H^\sigma(D;H)=\Nat j$ and $P_j$ is the $\Nat j$-orthogonal projection onto its range (Theorem~\ref{thm:kernel_LS}), the residual $e_j=(I-P_j)e_{j-1}$ satisfies $\|e_j\|_{\Nat j}\le\|e_{j-1}\|_{\Nat j}$. Combining this with \eqref{eq:scaled_a} and \eqref{eq:scaled_b} gives
\begin{equation}
\label{eq:step}
  \|e_j\|_{H^\sigma}\le\sqrt{c_2}\,\delta_j^{-\sigma}\|e_j\|_{\Nat j}
  \le\sqrt{c_2}\,\delta_j^{-\sigma}\|e_{j-1}\|_{\Nat j}
  \le C_2\,\delta_j^{-\sigma}\big(\|e_{j-1}\|_{L_2}+\delta_j^\sigma\|e_{j-1}\|_{H^\sigma}\big).
\end{equation}

\emph{Step 2: the $H^\sigma$ estimate for $j\ge2$.}
For $j\ge2$ the residual $e_{j-1}$ vanishes on $X_{j-1}$, so the sampling inequality \eqref{eq:SI_zero} gives $\|e_{j-1}\|_{L_2}\le C_Sh_{j-1}^\sigma\|e_{j-1}\|_{H^\sigma}$. Inserting this into \eqref{eq:step} and using $h_{j-1}/\delta_j=1/(\mu\nu)$ yields
\begin{equation}
\label{eq:contraction}
  \|e_j\|_{H^\sigma}\le C_2\Big(1+C_S\Big(\frac{h_{j-1}}{\delta_j}\Big)^{\sigma}\Big)\|e_{j-1}\|_{H^\sigma}
  =\beta\,\|e_{j-1}\|_{H^\sigma},\qquad j\ge2 .
\end{equation}

\emph{Step 3: the first level.}
For $j=1$ we have $e_0=B$, and \eqref{eq:step} with $\|B\|_{L_2}\le\|B\|_{H^\sigma}$ and $\delta_1\le1$ gives
\[
  \|e_1\|_{H^\sigma}\le C_2\,\delta_1^{-\sigma}\big(1+\delta_1^\sigma\big)\|B\|_{H^\sigma}\le2C_2\,\delta_1^{-\sigma}\|B\|_{H^\sigma}.
\]
Together with \eqref{eq:contraction},
\begin{equation}
\label{eq:eL_Hs}
  \|e_L\|_{H^\sigma}\le 2C_2\,\delta_1^{-\sigma}\,\beta^{\,L-1}\,\|B\|_{H^\sigma},\qquad L\ge1 .
\end{equation}

\emph{Step 4: proof of \eqref{eq:ms_L2} and \eqref{eq:ms_Linf}.}
The final residual $e_L=B-B_L$ vanishes on $X_L$, so \eqref{eq:SI_zero} gives
\[
  \|e_L\|_{L_2}\le C_S\,h_L^\sigma\,\|e_L\|_{H^\sigma}
  \qquad\text{and}\qquad
  \|e_L\|_{L_\infty}\le C_S\,h_L^{\sigma-d/2}\,\|e_L\|_{H^\sigma}.
\]
Now $h_L^\sigma\delta_1^{-\sigma}=(\mu^{L-1}h_1)^\sigma(\nu h_1)^{-\sigma}=\nu^{-\sigma}\mu^{\sigma(L-1)}$ and $\mu^{\sigma(L-1)}\beta^{L-1}=\alpha^{L-1}$, so \eqref{eq:eL_Hs} yields
\[
  \|e_L\|_{L_2}\le 2C_SC_2\,\nu^{-\sigma}\alpha^{L-1}\|B\|_{H^\sigma}
  \qquad\text{and}\qquad
  \|e_L\|_{L_\infty}\le 2C_SC_2\,\nu^{-\sigma}h_L^{-d/2}\alpha^{L-1}\|B\|_{H^\sigma},
\]
which are \eqref{eq:ms_L2} and \eqref{eq:ms_Linf}.

\emph{Step 5: the rates in (ii).}
If $\alpha<1$, then $\sigma_{\rm eff}=\log\alpha/\log\mu>0$, and $\sigma_{\rm eff}<\sigma$ because $\alpha>\mu^\sigma$, as $C_2\ge1$ and $C_S\nu^{-\sigma}>0$. Since $h_L/h_1=\mu^{L-1}$,
\[
  \alpha^{L-1}=\mu^{(L-1)\log\alpha/\log\mu}=\Big(\frac{h_L}{h_1}\Big)^{\sigma_{\rm eff}},
\]
and \eqref{eq:ms_L2} becomes the first bound in \eqref{eq:ms_order}. Likewise $h_L^{-d/2}\alpha^{L-1}=h_1^{-d/2}(h_L/h_1)^{\sigma_{\rm eff}-d/2}$ turns \eqref{eq:ms_Linf} into the second bound, and $\sigma_{\rm eff}>d/2$ is equivalent to $\alpha<\mu^{d/2}$.

\emph{Step 6: manifold-valued functions, part (iii).}
By Proposition~\ref{prop:transference}, $\|B-B_L\|_{L_2(D;H)}=\|F\ominus F_L\|_{L_2(D;G)}$ and $\|B-B_L\|_{L_\infty(D;H)}=\sup_x\dG(F(x),F_L(x))$. If $\alpha<\mu^{d/2}$, then $(h_L/h_1)^{\sigma_{\rm eff}-d/2}\le1$ for all $L$, so Step~5 gives $\sup_x\dG(F(x),F_L(x))\le\varepsilon:=2C_SC_2\nu^{-\sigma}h_1^{-d/2}\|B\|_{H^\sigma}$ uniformly in $L$, and Corollary~\ref{cor:transfer_natural} with this $\varepsilon$ transfers the $L_2$ rate in \eqref{eq:ms_order} to the natural distances and to the ambient Frobenius error.
\end{proof}

\section{Proofs of Propositions~\ref{prop:karcher} and~\ref{prop:karcher_general}}
\label{app:karcher}

\subsection{Proof of Proposition~\ref{prop:karcher}}

\begin{proof}
\emph{Step 1: strong convexity of $\phi$.}
Existence and uniqueness of $m_K$ in $B(e,2r)$ for small $r$, together with the characterisation
\[
  \grad\phi(m_K)=-\sum_iw_i\log_{m_K}(g_i)=0,
\]
are classical \cite{Karcher1977,Afsari2011}. For $g\in B(e,2r)$ one has $d(g,g_i)\le d(g,e)+d(e,g_i)\le2r+\|X_i\|\le3r$, using Theorem~\ref{thm:lie_metrics}(ii). Since the sectional curvature satisfies $K\le\kappa_+:=c_\g^2/4$, the Hessian comparison theorem gives, at distance $s$ from $q$,
\[
  \Hess\,\tfrac12d(\cdot,q)^2\;\ge\;\sqrt{\kappa_+}\,s\cot\big(\sqrt{\kappa_+}\,s\big)\,I ,
\]
and therefore $\phi$ is $c_0$-strongly geodesically convex on $B(e,2r)$ with
\[
  c_0:=3\sqrt{\kappa_+}\,r\cot\big(3\sqrt{\kappa_+}\,r\big)\ \ge\ \tfrac12\qquad\text{for } r\le r_0 ,
\]
with the continuous value $c_0=1$ when $\kappa_+=0$.

\emph{Step 2: from the gradient to the distance.}
For a $c_0$-strongly geodesically convex function $\phi$ with minimiser $m_K$ and any $g$ in the ball, strong convexity along the geodesic between $g$ and $m_K$ gives
\begin{align*}
  \phi(m_K)&\ge\phi(g)+\ip{\grad\phi(g)}{\log_gm_K}+\tfrac{c_0}{2}\,d(g,m_K)^2,\\
  \phi(g)&\ge\phi(m_K)+\tfrac{c_0}{2}\,d(g,m_K)^2 .
\end{align*}
Adding the two inequalities and using the Cauchy--Schwarz inequality yields
\begin{equation}
\label{eq:grad_to_dist}
  c_0\,d(g,m_K)\le\|\grad\phi(g)\| .
\end{equation}
Since $d(e,\bar g)\le\|\bar X\|\le r$ by Theorem~\ref{thm:lie_metrics}(ii), we may take $g=\bar g$.

\emph{Step 3: the gradient at the log-Euclidean mean.}
For a bi-invariant metric the geodesics through $\bar g$ are $t\mapsto\bar g\exp(tZ)$, so $\log_{\bar g}(g_i)=\mathrm dL_{\bar g}\log(\bar g^{-1}g_i)$ with the principal logarithm, and, as $\mathrm dL_{\bar g}$ is an isometry,
\[
  \|\grad\phi(\bar g)\|=\Big\|\sum_iw_i\log\big(e^{-\bar X}e^{X_i}\big)\Big\| .
\]
The map $(A,B)\mapsto\log(e^Ae^B)$ is analytic near $(0,0)$, and by the Baker--Campbell--Hausdorff formula
\[
  \log(e^Ae^B)=A+B+\tfrac12[A,B]+R(A,B),
  \qquad
  \|R(A,B)\|\le C_{\mathrm{BCH}}\big(\|A\|+\|B\|\big)^3 ,
\]
valid for $\|A\|+\|B\|\le r_1$, where $C_{\mathrm{BCH}}$ and $r_1$ depend only on $\g$ and its norm.
With $A=-\bar X$ and $B=X_i$, the linear terms sum to $\sum_iw_i(X_i-\bar X)=0$ and the quadratic terms to $-\frac12\big[\bar X,\sum_iw_iX_i\big]=-\frac12[\bar X,\bar X]=0$. Hence, for $2r\le r_1$,
\[
  \|\grad\phi(\bar g)\|\le\sum_iw_i\,\big\|R(-\bar X,X_i)\big\|\le C_{\mathrm{BCH}}\,(2r)^3 ,
\]
and \eqref{eq:grad_to_dist} with $c_0\ge\frac12$ gives $d(m_K,\bar g)\le16\,C_{\mathrm{BCH}}\,r^3$. Decreasing $r_0$ so that $2r_0\le r_1$ completes the proof.
\end{proof}

\subsection{Proof of Proposition~\ref{prop:karcher_general}}

\begin{proof}
\emph{Step 1: reduction to the identity.}
Let $g_0:=\exp X_0$ and let $\tau$ be the isometry of $(G,d)$ with $\tau(e)=g_0$: left translation $\tau(h)=g_0h$ for a bi-invariant metric, and the congruence $\tau(A)=g_0^{1/2}Ag_0^{1/2}$ for $\Spd_n$. Karcher means are equivariant under isometries, so $m_K=\tau(m_K')$, where $m_K'$ is the Karcher mean of the points $\tau^{-1}(g_i)=\exp Z(Y_i)$ with
\[
  Z(Y):=\log\big(\tau^{-1}(\exp(X_0+Y))\big),
\]
that is, $Z(Y)=\log(e^{-X_0}e^{X_0+Y})$ in the Lie-group case and $Z(Y)=\log(e^{-X_0/2}e^{X_0+Y}e^{-X_0/2})$ for $\Spd_n$.

\emph{Step 2: Taylor expansion of $Z$.}
The map $(X_0,Y)\mapsto Z(Y)$ is smooth on a neighbourhood of $\{\|X_0\|\le R\}\times\{0\}$ with $Z(0)=0$. Hence there are $r_1>0$ and constants $c_R$, $M_R$, uniform in $\|X_0\|\le R$, with
\[
  \|Z(Y)\|\le c_R\|Y\|\qquad\text{and}\qquad\|\mathrm d^2Z_Y\|\le M_R\qquad\text{for }\|Y\|\le r_1 .
\]
In particular $\|Z(Y_i)\|\le c_Rr$. Let $\bar Y:=\sum_iw_iY_i$ and $\bar Z:=\sum_iw_iZ(Y_i)$. A Taylor expansion of $Z$ at $\bar Y$, in which the linear terms cancel because the weights sum to one, gives
\begin{equation}
\label{eq:Z_taylor}
  \big\|\bar Z-Z(\bar Y)\big\|\le\frac{M_R}{2}\sum_iw_i\|Y_i-\bar Y\|^2\le2M_R\,r^2 .
\end{equation}

\emph{Step 3: the Karcher mean at the identity.}
Choose $r_1$ so small that $c_Rr_1\le r_0$, with $r_0$ from Proposition~\ref{prop:karcher}. That proposition, applied to the points $\exp Z(Y_i)$ with radius $c_Rr$, shows that $m_K'$ is the unique minimiser in $B(e,2c_Rr)$ and that
\[
  d\big(m_K',\exp\bar Z\big)\le C\,(c_Rr)^3 ;
\]
for $\Spd_n$ the same holds by the remark after Proposition~\ref{prop:karcher}.

\emph{Step 4: conclusion.}
Since $\bar g=\exp(X_0+\bar Y)=\tau(\exp Z(\bar Y))$ and $\tau$ is an isometry,
\begin{align*}
  d(m_K,\bar g)=d\big(m_K',\exp Z(\bar Y)\big)
  &\le d\big(m_K',\exp\bar Z\big)+d\big(\exp\bar Z,\exp Z(\bar Y)\big)\\
  &\le C\,c_R^3\,r^3+\Lam_0\,\big\|\bar Z-Z(\bar Y)\big\| .
\end{align*}
Here $\Lam_0=1$ by Theorem~\ref{thm:lie_metrics}(ii) for a bi-invariant metric, and $\Lam_0=\Lam(\sqrt2\,c_Rr_1)\le2$ by Theorem~\ref{thm:spd_metrics}(i) for $\Spd_n$, using $\omega(e^Z)\le\sqrt2\|Z\|_F$ for $\|Z\|_F\le c_Rr_1$ and decreasing $r_1$ if necessary. Together with \eqref{eq:Z_taylor},
\[
  d(m_K,\bar g)\le\big(C\,c_R^3\,r_1+2\Lam_0M_R\big)\,r^2 .\qedhere
\]
\end{proof}

\subsection{The example in Proposition~\ref{prop:karcher_general}}

\begin{proof}
Let $g_\pm:=\exp(X_0\pm tE)$ with $X_0=\diag(1,-1)$ and $E=E_{12}+E_{21}$, and let $J:=\diag(1,-1)$.

\emph{Step 1: the Karcher mean is diagonal.}
Since $JEJ=-E$ and $JX_0J=X_0$, conjugation by $J$ is an affine-invariant isometry that exchanges $g_+$ and $g_-$. The Karcher mean of two points with equal weights is the geodesic midpoint, here the geometric mean $M:=g_+\#g_-$, which is the unique positive definite solution of $Mg_+^{-1}M=g_-$. Being unique, $M$ is fixed by the conjugation, hence diagonal, and $\det M=(\det g_+\det g_-)^{1/2}=1$. Thus $M=\diag(e^s,e^{-s})$ for some $s\in\R$.

\emph{Step 2: an equation for $s$.}
With $\rho:=\sqrt{1+t^2}$ one has $(X_0\pm tE)^2=\rho^2I$, hence
\[
  g_\pm=\cosh\rho\,I+\frac{\sinh\rho}{\rho}\,(X_0\pm tE),\qquad
  g_+^{-1}=\cosh\rho\,I-\frac{\sinh\rho}{\rho}\,(X_0+tE).
\]
Comparing the $(1,1)$ entries in $Mg_+^{-1}M=g_-$ gives
\[
  e^{2s}\Big(\cosh\rho-\frac{\sinh\rho}{\rho}\Big)=\cosh\rho+\frac{\sinh\rho}{\rho},
  \qquad\text{that is,}\qquad
  s=\operatorname{artanh}\Big(\frac{\tanh\rho}{\rho}\Big)=:f(\rho).
\]

\emph{Step 3: expansion in $t$.}
Now $f(1)=1$ and
\[
  f'(1)=\big(\cosh^{-2}1-\tanh1\big)\cosh^21=1-\sinh1\cosh1=1-\tfrac12\sinh2 ,
\]
so with $\rho=1+\tfrac12t^2+\mathcal O(t^4)$ we obtain $s=1-\frac{\sinh2-2}{4}\,t^2+\mathcal O(t^4)$. Finally, since $M$ and $e^{X_0}$ are diagonal,
\[
  \dAI\big(M,e^{X_0}\big)=\big\|\log\big(e^{-X_0/2}Me^{-X_0/2}\big)\big\|_F=\sqrt2\,|s-1|
  =\frac{\sinh2-2}{2\sqrt2}\,t^2+\mathcal O(t^4),
\]
which is the claim.
\end{proof}

\end{document}